\documentclass[11pt,reqno]{amsart}
\usepackage{amsmath,amsfonts,amssymb,amscd,amsthm,amsbsy,epsf}
\newtheorem{thm}{Theorem} 
\newtheorem*{cor}{Corollary} 
\newtheorem*{CorThm}{Corollary to Theorem~\ref{groundstate}} 
\newtheorem{lem}[thm]{Lemma}
\newtheorem*{NSW}{Nonlinear Stein-Weiss Lemma}
\newtheorem*{NSW-groups}{Nonlinear Stein-Weiss Lemma on Groups}
\newtheorem*{SW-extended}{Stein-Weiss Lemma (extended)}
\newtheorem*{SymLem}{Symmetrization Lemma}
\newtheorem*{GSR}{Spectral representation for the fractional Laplacian 
(Frank-Lieb-Seiringer)}
\newtheorem*{HLS}{Hardy-Littlewood-Sobolev inequality}
\newtheorem*{HLS-b}{Hardy-Littlewood-Sobolev inequality}
\newtheorem*{loguncertain}{Logarithmic uncertainty}
\newtheorem*{logSobolev}{Logarithmic Sobolev inequality}
\newtheorem*{pitt}{Pitt's inequality}
\newtheorem*{pittspectral}{Pitt's inequality (at the spectral level)}
\newtheorem*{classical}{Classical Formula}
\newtheorem*{weighted}{Weighted Formula}
\newtheorem*{steinweiss}{Stein-Weiss inequality}

\theoremstyle{definition}

\newtheorem*{remarks}{Remarks}
\def\mysection#1{\par\medskip\noindent {\bf #1.}\par}
\def\ep{\varepsilon}
\def\F{{\mathcal F}}
\def\H{{\mathcal H}}
\def\S{{\mathcal S}}
\def\complex{{\mathbb C}}
\def\HH{{\mathbb H}}
\def\real{{\mathbb R}}
\def\Im{\operatorname{Im}}
\begin{document}

\title{Pitt's inequality and the fractional Laplacian:\\ sharp error estimates\\
$\quad$ \\ $\quad$ \\			%% to create 2 blank lines 
\hfill  $\text{\it for Eli Stein}$}
\author{William Beckner}
\address{Department of Mathematics, The University of Texas at Austin,
1 University Station C1200, Austin TX 78712-0257 USA}
\email{beckner@math.utexas.edu}
\begin{abstract}
Sharp error estimates in terms of the fractional Laplacian 
and a weighted Besov norm are obtained 
for Pitt's inequality by using the spectral representation with weights for the 
fractional Laplacian due to Frank, Lieb and Seiringer and the sharp 
Stein-Weiss inequality.
Dilation invariance, group symmetry on a non-unimodular group and a 
nonlinear Stein-Weiss lemma are used to provide short proofs of the 
Frank-Seiringer ``Hardy inequalities'' where fractional smoothness is 
measured by a Besov norm. 
\end{abstract}
%\dedicatory{for Eli Stein}
\maketitle

%%%%%%%%%%%%%%%%%%%%%%%%%%%%%
Considerable interest exists in understanding the framework of weighted 
inequalities for differential operators and the Fourier transform, and 
the application of quantitative information  drawn from these inequalities 
to varied problems in analysis and mathematical physics, including 
nonlinear partial differential equations,
spectral theory, fluid mechanics,
stability of matter, stellar dynamics  and uncertainty. 
Such inequalities provide both refined size estimates for differential 
operators and singular integrals and quantitative insight on symmetry 
invariance and geometric structure. 
The purpose of this note is to improve the sharp Pitt's inequality at the 
spectral level by using the optimal Stein-Weiss inequality and a new 
representation 
formula derived by Frank, Lieb and Seiringer \cite{FLS}
which expresses fractional Sobolev embedding in terms of a Besov norm 
characteristic of the problem's dilation invariance and extends with 
weights an earlier classical formula of Aronszajn and Smith using the 
$L^2$ modulus of continuity (see \cite{AS}, page~402).
Moreover reflecting the natural duality, 
this formula can be combined with the Hardy-Littlewood-Sobolev
inequality to provide new techniques to determine sharp embedding constants, 
including a sharp form of the Besov norm Sobolev embedding studied by 
Bourgain, Brezis and Mironescu \cite{BBM2000}. 
A secondary bootstrap argument produces an improved Stein-Weiss inequality 
accompanied by an intriguing monotonicity property at the spectral level. 
A direct proof of the weighted representation formula for the fractional 
Laplacian starting from the classical formula of Aronszajn and Smith 
is given in the appendix below. 
This classical formula provides an alternative way to view Pitt's inequality 
and Hardy inequalities in terms of fractional smoothness and Besov norms 
with natural extensions to $L^p$ classes. 
In this context, a new Stein-Weiss lemma, unusual in its simplicity, 
is applied to obtain the Frank-Seiringer ``Hardy inequalities'' for both 
$\real^n$ and the upper half-space $\real_+^n$ \cite{FS, FS-arXiv}, and 
corresponding extensions for the Heisenberg group and product spaces 
with mixed homogeneity.  

\mysection{1. Pitt's inequality}

This paper's principle result is the following theorem which gives 
specific improvement by a Besov norm for Pitt's inequality and explicitly 
demonstrates that the sharp constant is not attained which otherwise was 
observed from the equivalent Young's inequality (see \cite{B-PAMS95}).

\begin{thm}\label{groundstate}
For $f\in\S(\real^n)$ and $0<\beta <2$, $\beta \le \alpha <n$
\begin{equation}\label{eq:groundstate}
\begin{split}
C_\alpha \int_{\real^n} |\xi|^\alpha |\hat f|^2\,d\xi 
&\ge \int_{\real^n} |x|^{-\alpha} |f|^2 \,dx \\
\noalign{\vskip6pt}
&\qquad + \frac{C_\alpha}{D_\beta} \int_{\real^n\times\real^n} 
\frac{|g(x)- g(y)|^2}{|x-y|^{n+\beta}} 
|x|^{-(n-\beta)/2} |y|^{-(n-\beta)/2}\,dx\,dy
\end{split}
\end{equation}
where 
$$g(x) = |x|^{(n-\beta)/2} (-\Delta/4\pi^2)^{(\alpha-\beta)/4} f(x)$$
and 
\begin{gather*}
C_\alpha = \pi^\alpha 
\left[\frac{\Gamma (\frac{n-\alpha}4)}{\Gamma (\frac{n+\alpha}4)}\right]^2\ ,
\qquad 
D_\beta = \frac4{\beta} \ \pi^{\frac{n}2+ \beta} \ 
\frac{\Gamma (1-\frac{\beta}2)}{\Gamma (\frac{n+\beta}2)}\ .
%\\
%\noalign{\vskip6pt}
%D_{\alpha,\beta} = \pi^{-n/2\  +\alpha -\beta}\ 
%\frac{\beta}4\ \frac{\Gamma (\frac{n+\beta}2)}{\Gamma (\frac{2-\beta}2)}
%\left[\frac{\Gamma (\frac{n-\alpha}4)}{\Gamma (\frac{n+\alpha}4)}
%\right]^2\ .
\end{gather*}
\end{thm}

Before outlining the proof for this result, the context for the problem will
be described. 
Note that equation~\eqref{eq:groundstate} becomes an equality if 
$\alpha=\beta$, the Fourier transform is defined by 
$$(\F f) (\xi) = \hat f(\xi) = \int_{\real^n} e^{2\pi i\xi x} f(x)\,dx\ ,$$
primes denote dual exponents, $1/p + 1/p' =1$, and $\Delta$ is the 
Laplacian on $\real^n$.

The framework for this development lies at the interface of two 
classical results combined with the recent representation formula for 
the fractional Laplacian obtained by Frank, Lieb and Seiringer: 
namely,  Pitt's inequality which may be viewed as a weighted extension of 
the Hausdorff-Young inequality and the Stein-Weiss inequality which 
extends the Hardy-Littlewood-Sobolev inequality for fractional integrals 
and in effect Young's inequality for convolution 
thus reinforcing the natural operator duality.
This paper is the fourth in a series \cite{B-PAMS95, B-weighted06, B-pitt06} 
which was originally motivated by an effort to better understand how the 
classical uncertainty principle could be extended to incorporate the 
logarithmic Sobolev inequality. 
In a larger sense, the purpose is to determine the implications of dilation 
invariance on the Euclidean manifold. 

\begin{pitt}
For $f\in \S(\real^n)$, $1<p\le q,\infty$, $0\le \alpha <n/q$, 
$0\le \beta < n/p'$ and $n\ge2$
\begin{equation}\label{eq:pitt}
\bigg[ \int_{\real^n} \big|\, |x|^{-\alpha} \hat f\big|^q\,dx\bigg]^{1/q} 
\le A\bigg[ \int_{\real^n} \big|\, |x|^\beta f\big|^p\,dx\bigg]^{1/p}
\end{equation}
with the index constraint 
$$\frac{n}p + \frac{n}q +\beta -\alpha =n\ .$$
\end{pitt}

\begin{pittspectral}
For $f\in \S(\real^n)$ and $0\le \alpha <n$
\begin{gather}
\int_{\real^n} |x|^{-\alpha} |f(x)|^2\,dx 
\le C_\alpha \int_{\real^n} |\xi|^\alpha |\hat f(\xi)|^2\,d\xi
\label{eq:pittspectral}\\
\noalign{\vskip6pt}
C_\alpha = \pi^\alpha \left[ \Gamma \left(\frac{n-\alpha}4\right) \Big/
\Gamma \left(\frac{n+\alpha}4\right)\right]^2 \notag
\end{gather}
\end{pittspectral}

\begin{steinweiss}\label{thm:steinweiss}
For $f\in L^p(\real^n)$ with $1<p<\infty$, $0<\lambda <n$, $\alpha <n/p$, 
$\beta <n/p'$ and $n=\lambda +\alpha +\beta$
\begin{gather}
\left\||x|^{-\alpha} (|x|^{-\lambda} * (|x|^{-\beta} f))\right\|_{L^p(\real^n)}
\le E_{\alpha,\beta} \|f\|_{L^p (\real^n)}
\label{eq:steinweiss}\\
\noalign{\vskip6pt}
E_{\alpha,\beta} = \pi^{n/2} 
\left[\frac{\Gamma (\frac{\alpha+\beta}2) 
\Gamma (\frac{n}{2p} - \frac{\alpha}2) 
\Gamma (\frac{n}{2p'} - \frac{\beta}2)}
{\Gamma (\frac{n-\alpha-\beta}2) 
\Gamma (\frac{n}{2p'} + \frac{\alpha}2) 
\Gamma (\frac{n}{2p} + \frac{\beta}2)} \right] 		\notag
\end{gather}
\end{steinweiss}

\noindent
Note that it is not required that both $\alpha$ and $\beta$ are 
non-negative, but rather $0<\alpha +\beta <n$.

\begin{GSR}
For \newline $f\in \S(\real^n)$ and $0\le \alpha <\min \{2,n\}$ 
\begin{gather}
C_\alpha \int_{\real^n} |\xi|^\alpha |\hat f(\xi)|^2\,d\xi
= \int_{\real^n} |x|^{-\alpha} |f(x)|^2\,dx 
\label{eq:GSR}\\
\noalign{\vskip6pt}
+ D_\alpha \int_{\real^n\times\real^n} 
\frac{|g(x) - g(y)|^2}{|x-y|^{n+\alpha}} 
|x|^{-(n-\alpha)/2} |y|^{-(n-\alpha)/2}\, dx\,dy		\notag
\end{gather}
where 
\begin{gather*}
g(x) = |x|^{(n-\alpha)/2} f(x)\ ,\\
\noalign{\vskip6pt}
C_\alpha = \pi^\alpha \left[ \Gamma \left(\frac{n-\alpha}4\right) \Big/
\Gamma \left( \frac{n+\alpha}4\right)\right]^2\ ,\\
\noalign{\vskip6pt}
D_\alpha = \pi^{-n/2\, -\,\alpha}\  
\frac{\alpha}4\  \frac{\Gamma (\frac{n+\alpha}2)}{\Gamma (\frac{2-\alpha}2)}
\ C_\alpha\ .
\end{gather*}
\end{GSR}

\noindent
Observe that in the limit $\alpha\to2$, one obtains the classical 
relation which includes Hardy's inequality 
\begin{equation}\label{classicalinequality}
\int_{\real^n} |\nabla f|^2\,dx 
= \frac{(n-2)^2}4 \int_{\real^n} |x|^{-2} |f(x)|^2\,dx 
+ \int_{\real^n} |\nabla (|x|^{(n-2)/2} f)|^2 |x|^{-n+2} \,dx
\end{equation}

\begin{proof}[Proof of Theorem~\ref{groundstate}]
Use the ground state representation \eqref{eq:GSR} to write for $\alpha >\beta$
\begin{gather*}
C_\beta \int_{\real^n} |\xi|^\beta \left[ |\xi|^{\alpha-\beta} |\hat f|^2
\right]\,d\xi 
=\int_{\real^n}|x|^{-\beta} |\F^{-1} (|\xi|^{(\alpha-\beta)/2}\hat f\,)|^2\,dx\\
\noalign{\vskip6pt}
+ D_\beta \int_{\real^n\times\real^n} 
\frac{|g(x) -g(y)|^2}{|x-y|^{n+\beta}} \ 
|x|^{-(n-\beta)/2} |y|^{-(n-\beta)/2}\,dx\,dy
\end{gather*}
where 
$$g(x) = |x|^{(n-\beta)/2} (-\Delta/4\pi^2)^{(\alpha-\beta)/4} f(x)\ .$$
To complete the proof, it remains to show that 
\begin{equation}\label{eq:groundstateproof}
\frac{C_\beta}{C_\alpha} \int |x|^{-\alpha} |f|^2\, dx
\le\int_{\real^n}|x|^{-\beta}|\F^{-1}
(|\xi|^{(\alpha-\beta)/2}\hat f\,)|^2\,dx\ .
\end{equation}
But this is equivalent to showing the Stein-Weiss inequality 
\begin{equation*}
\int_{\real^n} |x|^{-\alpha} \left| \, |x|^{-n+(\alpha-\beta)/2} * 
(|x|^{\beta/2} f)\right|^2\, dx 
\le F_{\alpha,\beta} \int_{\real^n} |f|^2\,dx
\end{equation*}
with 
\begin{gather*}
F_{\alpha,\beta} = \frac{C_\alpha}{C_\beta}\ \pi^{n-\alpha+\beta} 
\left[ \frac{\Gamma (\frac{\alpha-\beta}4)}{\Gamma (\frac{2n-\alpha+\beta}4)}
\right]^2\\
\noalign{\vskip6pt}
= \pi^n \left[ \frac{\Gamma (\frac{\alpha-\beta}4) 
\Gamma (\frac{n-\alpha}4) \Gamma (\frac{n+\beta}4)} 
{\Gamma (\frac{2n-\alpha+\beta}4) 
\Gamma (\frac{n+\alpha}4) \Gamma (\frac{n-\beta}4)} \right]^2\ .
\end{gather*}
This constant exactly matches the bound given by the Stein-Weiss inequality 
\eqref{eq:steinweiss} though note the change in notation and that here it is a 
positive power of $|x|$ multiplying the function $f$.
Then 
\begin{equation*}
D_{\alpha,\beta} = C_\alpha D_\beta/C_\beta  
= \pi^{-n/2\ +\alpha-\beta} \ \frac{\beta}4\ 
\frac{\Gamma (\frac{n+\beta}2)}{\Gamma (\frac{2-\beta}2)}
\left[\frac{\Gamma (\frac{n-\alpha}4)}{\Gamma (\frac{n+\alpha}4)}\right]^2\ .
\end{equation*}
This completes the proof of Theorem~\ref{groundstate}.
\renewcommand{\qed}{}
\end{proof}

As evident from the proof outlined above and the previous determination of 
sharp constants for Pitt's inequality, there is a self-evident duality 
between Pitt's inequality and the Stein-Weiss inequality. 
Using 
$$\F [|x|^{-\alpha/2}] = \sqrt{d_\alpha} \ |x|^{-n+\frac{\alpha}2}\ ,
\qquad 
d_\alpha = \pi^{-(n-\alpha)} 
\left[ \frac{\Gamma (\frac{2n-\alpha}4)}{\Gamma (\frac{\alpha}4)}\right]^2$$
then the following error estimate for the Stein-Weiss inequality is 
obtained from Theorem~\ref{groundstate}: 

\begin{thm}\label{error-est:stein-weiss} 
For $f\in L^2 (\real^n)$ and $0<\beta <2$, $\beta \le \alpha <n$ 
\begin{equation}\label{eq:stein-weiss} 
\begin{split}
\int_{\real^N} |f|^2\,dx 
&\ge \frac{d_\alpha}{C_\alpha} \int_{\real^n} |x|^{-\alpha} 
\Big| \, |x|^{-n+\frac{\alpha}2} * f\Big|^2\,dx \\
\noalign{\vskip6pt}
&\qquad + \frac{d_\beta}{D_\beta} \int_{\real^n\times \real^n} 
\frac{|h(x) - h(y)|^2}{|x-y|^{n+\beta}} \ 
|x|^{-(n-\beta)/2} \ |y|^{-(n-\beta)/2}\, dx\,dy\ .
\end{split}
\end{equation}
where
$$h(x) = \left(  |x|^{-n+\frac{\beta}2} * f\right) (x)\ .$$ 
\end{thm}

\noindent
The expression explicitly demonstrates that the sharp constant for the 
Stein-Weiss inequality is not attained. 
Using  the Frank-Lieb-Seiringer representation formula, one obtains the 
following monotonicity result:

\begin{cor}
For $0<\alpha <2$ 
$$\frac{d_\alpha}{C_\alpha} \int_{\real^n} |x|^{-\alpha} \Big|\, 
|x|^{-n+\frac{\alpha}2} * f\Big|^2\,dx$$
is monotone decreasing as a function of $\alpha$!
\end{cor}

A natural objective in transforming weighted inequalities for the Fourier
transform to convolution estimates is to enable extensions to 
$L^p (\real^n)$ (see \cite{B-pitt06} where this aspect was the central 
direction). 
The classical Aronszajn-Smith representation formula (see Appendix in this 
paper) allows Theorem~\ref{groundstate} to be recast in a form that is 
amenable to analysis on $L^p(\real^n)$ (see section~4 below). 

\begin{thm}\label{thm:recast}
For $f\in \S(\real^n)$ and $0<\beta <2$, $\beta \le\alpha <n$ 
\begin{equation}\label{eq:recast}
\begin{split}
&\frac{C_\alpha}{D_\beta} \int_{\real^n\times\real^n} 
\frac{|h(x) - h(y)|^2}{|x-y|^{n+\beta}} \,dx\,dy 
\ge \int_{\real^n} |x|^{-\alpha} |f|^2\,dx\\
\noalign{\vskip6pt}
&\qquad 
+ \frac{C_\alpha}{D_\beta} \int_{\real^n\times\real^n} 
\frac{|g(x) - g(y)|^2}{|x-y|^{n+\beta}} \ 
|x|^{-(n-\beta)/2}\ |y|^{-(n-\beta)/2}\,dx\,dy 
\end{split}
\end{equation}
with $h(x) = |x|^{-(n-\beta)/2} g(x) 
= (-\Delta/4\pi^2)^{(\alpha-\beta)/4}f(x)$.
\end{thm}

\mysection{2. Logarithmic uncertainty}
\medskip

The central motivation for computing optimal constants for Pitt's inequality 
in \cite{B-PAMS95} 
was to directly strengthen the classical uncertainty principle.
Beginning with analysis of the log Sobolev inequality, one understands 
that a parameter-dependent variational inequality which becomes
an equality at a given parameter value can be differentiated at that value 
to give a new analytic inequality. 
Using Pitt's inequality, one obtains 
(see \cite{B-PAMS95}, \cite{B-pitt06}): 

\begin{loguncertain}
For $f\in \S(\real^n)$ and $1<p<\infty$ 
\begin{gather}
\int_{\real^n}\ln |x|\,|f|^2\,dx + \int_{\real^n} \ln |\xi|\, |\hat f|^2\,d\xi 
\ge D\int_{\real^n} |f|^2\,dx 
\label{loguncertain1}\\
\noalign{\vskip6pt}
D = \psi (n/4) - \ln \pi\ ,\qquad \psi = (\ln \Gamma )' \notag\\
\noalign{\vskip6pt}
\int_{\real^n} \ln |x|\, |f|^p\, dx 
+\int_{\real^n}\left[(\ln \sqrt{(-\Delta/4\pi^2)}\, )f\right] f |f|^{p-2}\, dx
\ge E\int_{\real^n} |f|^p\,dx
\label{loguncertain2}\\
\noalign{\vskip6pt}
E = \frac12 \Big[ \psi (n/2p) + \psi(n/2p')\Big] - \ln \pi \notag
\end{gather}
\end{loguncertain}

\noindent
Since the Frank-Lieb-Seiringer spectral formula is a direct relation, it can 
be differentiated for all allowed values of $\alpha$, including $\alpha=0$.
This argument demonstrates that no extremals exist for logarithmic 
uncertainty.

\begin{thm}\label{thm2}
For $f\in \S (\real^n)$ and $g= |x|^{n/2} f$
\begin{equation} \label{eq8}
\begin{split}
&\int_{\real^n} \ln |x|\, |f|^2\,dx 
+ \int_{\real^n} \ln |\xi |\, |\hat f|^2\,d\xi  - D\int_{\real^n} |f|^2\,dx \\
\noalign{\vskip6pt}
&\qquad = \frac14 \pi^{-n/2}\ \Gamma (n/2) \int_{\real^n\times\real^n} 
\frac{|g(x) - g(y)|^2}{|x-y|^n}\ |x|^{-n/2} |y|^{-n/2}\,dx\,dy\ .
\end{split}
\end{equation}
with $D = \psi (n/4) - \ln \pi$.
\end{thm}

This result can be compared to the conformally invariant 
Hardy-Littlewood-Sobolev inequality for $L^2 (\real^n)$ 
(see Theorem~3 in \cite{B-PAMS95}):

\begin{HLS}[$L^2$ entropy form]\label{thm:HSL}
For $f\in \S(\real^n)$ with $\|f\|_2=1$ 
\begin{gather}
\frac{n}2 \int_{\real^n} \ln |\xi|\, |\hat f(\xi)|^2\,d\xi 
\ge \int_{\real^n} \ln |f(x)|\, |f(x)|^2\,dx + B_n 
\label{eq:HLS}\\
\noalign{\vskip6pt}
B_n = \frac{n}2 \ \psi \left(\frac{n}2\right) - \frac12 \ln 
\left[\pi^{n/2} \ \Gamma (n)/\Gamma (n/2)\right] \notag
\end{gather}
Up to conformal automorphism, extremal functions are of the form 
$A(1+|x|^2)^{-n/2}$.
\end{HLS}

But observe that a weaker form of this inequality is given by ``logarithmic 
uncertainty'' \eqref{loguncertain1} which will suffice by an asymptotic 
argument to determine the logarithmic Sobolev inequality in its Euclidean 
form (see \cite{B-ForumMath99}, page~117).

\begin{logSobolev}			%% {Logarithmic Sobolev inequality}
For $f\in \S(\real^n)$ with $\|f\|_2 =1$
\begin{equation}\label{eq:logSobolev}
\int_{\real^n} |f|^2 \ln |f|\,dx 
\le \frac{n}4 \ \ln \bigg[ \frac2{\pi en} \int_{\real^n} |\nabla f|^2\,dx 
\bigg]\ .
\end{equation}
\end{logSobolev}

\begin{thm}\label{thm:logSobolev}
Logarithmic uncertainly implies the logarithmic Sobolev inequality.
\end{thm}

\begin{proof}
Observe that for a non-negative radial decreasing function $h(x)$ 
with $\|h\|_2 =1$, $h(x) \le m(B)^{-1/2} |x|^{-n/2}$ where $m(B)$ is 
the volume of the unit ball. 
Then 
$$\int_{\real^n} \ln |x|\, |h|^2\,dx 
\le -\frac2n \int_{\real^n} \ln |h|\, |h|^2\,dx 
- \frac1n \ln \left[ \pi^{n/2}/\Gamma \Big(\frac{n}2 +1\Big)\right]$$
and in general
$$\int_{\real^n} \ln |\xi|\, |\hat h|^2\,d\xi 
\le \frac12 \ln \int_{\real^n} |\nabla h|^2\,dx - \ln (2\pi)$$
so that 
\begin{gather*}
\frac12 \ln \int_{\real^n} |\nabla h|^2\, dx 
- \frac2n \int_{\real^n} \ln |h|\, |h|^2\, dx 
-\ln(2\pi) - \frac1n \ln \left[ \pi^{n/2}/\Gamma\Big(\frac{n}2 +1\Big)\right]\\
\noalign{\vskip6pt}
\ge \int_{\real^n} \ln |\xi|\, |\hat h|^2\, d\xi 
+ \int_{\real^n} \ln |x|\, |h|^2\,dx 
\ge \psi(n/4) - \ln \pi
\end{gather*}
and 
$$\ln \int_{\real^n} |\nabla h|^2\,dx -\frac4n \int \ln |h|\, |h|^2\,dx 
\ge 2\psi (n/4) +\ln (4\pi) - \frac2n \ln \, \Gamma \Big(\frac{n}2 +1\Big)\ .$$
Since this result will hold for the equimeasurable radial decreasing 
rearrangement of a function and the left-hand side of this inequality decreases
under symmetric rearrangement, it must hold for all functions $F$ with 
$\|F\|_2 =1$:
\begin{gather}
\ln \left[ \frac2{\pi en} \int_{\real^n} |\nabla F|^2\,dx \right]
- \frac4n \int_{\real^n} \ln |F|\, |F|^2\,dx \label{eq:logSobolev1}\\
\noalign{\vskip6pt}
\ge 2\psi (n/4) + \ln (8/en) - \frac2n \ln \, \Gamma \Big(\frac{n}2+1\Big)\ .
\notag
\end{gather}
Set $n=km$ with $F(x) = \prod f(x_j)$ for $j=1$ to $k$ with $x_j\in\real^m$ 
and $\|f\|_2=1$. 
Then \eqref{eq:logSobolev1} becomes 
\begin{equation*}
\begin{split}
&\ln \left[ \frac2{\pi em} \int_{\real^m} |\nabla f|^2\,dx \right]
- \frac4m \int_{\real^m} \ln |f|\, |f|^2\,dx\\
\noalign{\vskip6pt}
&\ge 2\psi \left( \frac{mk}4\right) + \ln \left[\frac8{emk}\right] 
- \frac2{mk} \ln\, \Gamma \left(\frac{mk}2 +1\right)\\
\noalign{\vskip6pt}
&\qquad \simeq -\frac1{mk} \ln (\pi mk) \to 0
\end{split}
\end{equation*}
as $k\to0$. 
Hence 
$$\ln \left[ \frac2{\pi em} \int_{\real^m} |\nabla f|^2\,dx \right] 
- \frac4m \int_{\real^m} \ln |f|\, |f|^2\,dx \ge 0$$
and the proof of Theorem~\ref{thm:logSobolev} is complete.
A brief sketch of this argument was given in \cite{B-ForumMath99}.
\renewcommand{\qed}{}
\end{proof}

By utilizing the condition for equality in Theorem~\ref{groundstate} 
$(\alpha =\beta)$, a second extension of logarithmic uncertainty 
is obtained through a differentiation argument.

\begin{CorThm}
For $f\in \S(\real^n)$ and $0\le\beta <2$ 
\begin{equation} \label{eq:CorThm}
\begin{split}
&C_\beta \int_{\real^n} |\xi|^\beta \ln |\xi|\, |\hat f|^2\,d\xi 
+ \int_{\real^n} |x|^{-\beta} \ln |x|\, |f|^2\,dx\\
\noalign{\vskip6pt}
&\hskip.5in \ge \left( \frac12 \Big[ \psi \Big( \frac{n+\beta}4\Big) 
+ \psi \Big(\frac{n-\beta}4\Big)\Big] - \ln \pi\right) 
\int_{\real^n} |x|^{-\beta} |f|^2\,dx \ .
\end{split}
\end{equation}
\end{CorThm}

\mysection{3. Besov norms and the Hardy-Littlewood-Sobolev inequality}
\medskip

The nature of the weighted Besov norm can best be understood in terms 
of dilation invariance which allows its reformulation on the 
multiplicative group $\real_+$ or the real line $\real$. 
Combine equation~\eqref{eq8} with the $L^2$ Hardy-Littlewood-Sobolev 
inequality to obtain for $\|f\|_2 =1$
\begin{gather}
\frac14 \pi^{-n/2} \Gamma (n/2)\int_{\real^n\times\real^n} 
\frac{|g(x)-g(y)|^2}{|x-y|^n} \ |x|^{-n/2} |y|^{-n/2}\,dx\,dy
\label{weightedBesov}\\
\noalign{\vskip6pt}
\ge \frac2n \int_{\real^n}  \ln |g|\, |f|^2\, dx + E_n\notag\\
\noalign{\vskip6pt}
E_n = \psi (n/2) - \psi (n/4) - \frac1n  \ln 
\left[\pi^{-n/2} \Gamma (n) \Big/ \Gamma (n/2)\right]\ .\notag
\end{gather}
Observe that due to equality in equation~\eqref{eq8}, 
inequality~\eqref{weightedBesov} is improved when $f$ is replaced by 
its equimeasurable radial decreasing rearrangement. 
Hence $f$ can be taken as radial and set $t= |x|$, $g(t) = g(x)$. 
Then
\begin{gather}
\frac14 \pi^{-n/2} \Gamma (n/2) \int_{\real_+\times\real_+} 
|g(s) - g(t)|^2\, \varphi \Big[ \frac{s}t +\frac{t}s -2\Big] 
\, \frac{ds}s \, \frac{dt}t
\label{f-radial}\\
\noalign{\vskip6pt}
\ge \frac2n \int_{\real_+} \ln |g|\, |g|^2\ \frac{dt}t + E_n\notag
\end{gather}
and

\begin{thm}\label{thm4}
For $v\in \S(\real)$ with $\|v\|_2 =1$
\begin{gather}
\frac14 \pi^{-n/2} \Gamma (n/2) \int_{\real\times\real} 
|v(x) - v(y)|^2\ \varphi\! \left[4\sinh^2 \Big(\frac{x-y}2\Big)\right] \,dx\,dy
\label{eq:thm4}\\
\noalign{\vskip6pt}
\ge \frac2n \int_{\real} \ln |v|\, |v|^2 \,dx + E_n\notag\\
\noalign{\vskip6pt}
\varphi (w) = \int_{S^{n-1}} [w+2(1-\xi_1)]^{-n/2}\,d\xi \notag
\end{gather}
where $d\xi$ is normalized surface measure.
Extremal functions include $v(x) = A(\cosh x)^{-n/2}$.
\end{thm}

The full Hardy-Littlewood-Sobolev inequality can be used with the 
Frank-Lieb-Seiringer formula to give an embedding result on the 
multiplicative group in terms of the Besov norm.

\begin{HLS-b}
For $f\in L^p (\real^n)$, $1<p<2$
\begin{gather}
\Big| \int_{\real^n\times\real^n} f(x) |x-y|^{-2n/p'} g(y)\,dx\,dy\Big|
\le A_p \|f\|_{L^p(\real^n)} \|g\|_{L^p (\real^n)}\label{eq:HLS-b1}\\
\noalign{\vskip6pt}
A_p  = \pi^{n/p'}\ \frac{\Gamma (\frac{n}p - \frac{n}2)}{\Gamma (\frac{n}p)}
\ \left[ \frac{\Gamma (n)}{\Gamma (\frac{n}2)}\right]^{2/p\, -\, 1} \notag\\
\noalign{\vskip6pt}
[ \|f\|_{L^{p'} (\real^n)} ]^2 
\le C_p \int_{\real^n} |(-\Delta)^{n(2/p\, - 1)/4} f|^2\,dx 
= (2\pi)^{n(2/p\, -\,1)} C_p \int_{\real^n} |\xi|^{n(2/p\,-\,1)} |\hat f|^2\, 
d\xi \label{eq:HLS-b2}\\
\noalign{\vskip6pt}
C_p  = \frac{\Gamma (\frac{n}{p'})}{\Gamma (\frac{n}p)}
\ \left[\frac{\Gamma (n)}{(4\pi)^{n/2}\Gamma (\frac{n}2)}\right]^{2/p\, -\,1}
\ . \notag
\end{gather}
\end{HLS-b}

Set $\alpha = n(2/p -1)$; then from \eqref{eq:GSR} for $0<\alpha<\min\{2,n\}$
and $g = |x|^{(n-\alpha)/2} f$ 
\begin{gather}
C_\alpha \frac{\Gamma (\frac{n+\alpha}2)}{\Gamma (\frac{n-\alpha}2)} \ 
\left[ \frac{\Gamma (\frac{n}2)}{\Gamma (n)}\right]^{\alpha/n} \ 
[\|f\|_{L^{p'}  (\real^n)} ]^2 
\le \int_{\real^n} |x|^{-\alpha} |f|^2\,dx \label{eq:GSR-b} \\
\noalign{\vskip6pt}
+ D_\alpha \int_{\real^n\times\real^n} 
\frac{|g(x) - g(y)|^2}{|x-y|^{n+\alpha}} \
|x|^{-(n-\alpha)/2} |y|^{-(n-\alpha)/2}\,dx\,dy\ . \notag
\end{gather}
Note that for $n=2$ the condition on $\alpha$ allows the full range of 
values for $1<p<2$. 
Since $f$ can be taken to be radial, let $t= |x|$ and 
$f(x) = g(t) t^{-(n-\alpha)/2}$; then one obtains the following Besov estimate
on $\real_+$ for $2<p' <p_c$ where  $p_c = 2n/(n-2)$ is the critical 
Sobolev embedding index on $\real^n$ for $n>2$.

\begin{thm}\label{thm:Besov-estimate}
For $g\in L^{p'}(\real_+)$ with $\alpha = n(1-2/p')$ and $n>2$, 
$2<p' < 2n/(n-2)$
\begin{equation}\label{eq:Besov-estimate}
\begin{split}
&C_\alpha\ \frac{\Gamma (\frac{n+\alpha}2)}{\Gamma (\frac{n-\alpha}2)}
\left[ \frac{\Gamma (\frac{n}2) \Gamma (\frac{n}2)} {2\pi^{n/2} \Gamma (n)}
\right]^{\alpha/n} \|g\|_{L^{p'}(\real_+)} \\
\noalign{\vskip6pt}
&\qquad \le \Big[ \|g\|_{L^2 (\real_+)}\Big]^2 
+ D_\alpha \left[ \frac{2\pi^{n/2}}{\Gamma (\frac{n}2)}\right] 
\int_{\real_+\times \real_+}  \mkern-24mu
|g(t) - g(s)|^2 \psi_{\alpha,n} (t/s)\ \frac{dt}t\ \frac{ds}s
\end{split}
\end{equation}
where 
$$\psi_{\alpha,n} (t) 
= \int_{S^{n-1}} |t+\frac1t -2\xi_1|^{-(n+\alpha)/2} \, d\xi$$
and 
\begin{gather*}
C_\alpha = \pi^\alpha \left[ \Gamma \Big(\frac{n-\alpha}4\Big) \Big/ 
\Gamma \Big(\frac{n+\alpha}4\Big)\right]^2\ ,\\
\noalign{\vskip6pt}
D_\alpha = \pi^{-n/2 -\alpha} \ \frac{\alpha}4 \ 
\frac{\Gamma (\frac{n+\alpha}2)}{\Gamma (\frac{2-\alpha}2)} \ C_\alpha\ .
\end{gather*}
\end{thm}
\medskip

The classical Aronszajn-Smith representation formula for the fractional 
Laplacian combined with the Hardy-Littlewood-Sobolev inequality allows
one to obtain the sharp $L^2$ embedding constant for the Besov norm 
inequality considered by Bourgain, Brezis and Mironescu 
\cite{BBM2000, BBM02}, and Maz'ya and Shaposhnikova \cite{MS}.

\begin{thm}\label{thm:dualform}
For $F\in \S(\real^n)$, $0<\beta <1$, $2<q = 2n/(n-2\beta) < p_c= 2n/(n-2)$
\begin{align} 
\left[\|f\|_{L^q(\real^n)}\right]^2 
& \le C(n,2,\beta) \int_{\real^n\times\real^n} 
\frac{|f(x) - f(y)|^2}{|x-y|^{n+2\beta}}\,dx\,dy \label{eq:dualform1}\\
\noalign{\vskip6pt}
C(n,2,\beta) & = \frac{\beta (1-\beta)}{n-2\beta} \ \pi^{-\beta -n/2}
\ \frac{\Gamma (\frac{n}2 +1-\beta)}{\Gamma(2-\beta)} 
\ \left[\frac{\Gamma (n)}{\Gamma (\frac{n}2)} \right]^{2\beta/n}\ .
\notag
\end{align}
\end{thm}

\begin{proof} 
This result follows directly from the dual form of the Hardy-Littlewood-Sobolev
inequality for the fractional Laplacian given by equation~\ref{eq:HLS-b2} 
(also see equation~(18) in reference \cite{B-ForumMath99}) 
and the classical Aronszajn-Smith representation formula from the Appendix
below (see \eqref{eq:weighted}) 
by setting $p'=q$ and $2\beta = n(2/p -1)$. 
\end{proof}

This inequality is sharp since for the allowed range of parameters it is 
equivalent to the Hardy-Littlewood-Sobolev inequality and one can 
calculate extremal functions using the extremals for that inequality.
Observe that $C(n,2,\beta)$ has the appropriate character determined 
by the respective principal theorems in Bourgain, Brezis and Mironescu
(see equation (8), page 78 in \cite{BBM2000}) and in 
Maz'ya and Shaposhnikova (see equation (2), page 231 in \cite{MS}). 

%%%%%%%%%%%%%%%%%%%%%%%%%%
\mysection{4. Besov norms and the nonlinear Stein-Weiss lemma}

In the study and understanding of functional inequalities, symmetry often
plays a central role, particularly in identifying optimal constants and
extremal functions.
A useful conceptual tool has been to transfer symmetry structures to analysis
on a group (see \cite{B-Geo} and \cite{B-PAMS01}). 
The classical  example is the Hardy-Littlewood paradigm that a positive
integral operator that commutes with dilations may be reduced to Young's
inequality for sharp convolution estimates on the multiplicative group.
This principle is demonstrated here by giving  short elementary proofs 
of two recent results by Frank and Seiringer 
\cite{FS, FS-arXiv} on weighted inequalities
and Besov norms to measure fractional smoothness.
Symmetry invariance allows reduction of the problem to convolution 
estimates on a Lie group where a nonlinear Stein-Weiss lemma is used 
to obtain optimal estimates. 

\begin{thm}\label{thm8}
Let $f\in \S(\real^n)$, $0<\beta <1$ and  $1\le p< n/\beta$; then
\begin{gather}
\int_{\real^n\times \real^n}
\frac{|f(x) - f(y)|^p}{|x-y|^{n+p\beta}}\,dx\,dy
\ge D_{p,\beta} \int_{\real^n} |x|^{-p\beta} |f(x)|^p\,dx\label{eq1}\\
\noalign{\vskip6pt}
D_{p,\beta} = \int_{\real^n} \big|1-|x|^{-\lambda}\big|^p
|x-\eta|^{-n-p\beta}\,dx
\notag
\end{gather}
for $\lambda = (n-p\beta)/p$ and $\eta \in S^{n-1}$.
\end{thm}

The proof depends on using the dilation invariance to transfer this
inequality to a convolution problem on the multiplicative group and
application of a new nonlinear Stein-Weiss lemma (see the appendix in
\cite{B-pitt06} for the classical Stein-Weiss lemma).

\begin{NSW}
For $f,g\in L^p(\real^m)$, $1\le p <\infty$
\begin{align}
&\int_{\real^m\times\real^m} |g(y-x)f(x) - g(x-y)f(y)|^p\, dx\,dy
\notag\\
\noalign{\vskip6pt}
&\qquad
\ge \int_{\real^m} \big|\,|g(y)| - |g(-y)|\,\big|^p\,dy
\int_{\real^m} |f(x)|^p\,dx\ .
\label{eq:NSW}
\end{align}
\end{NSW}

\begin{proof}
By a change of variables in $y$ and using the triangle inequality for norms:
\begin{align*}
&\int_{\real^m\times\real^m} |g(y-x)f(x) - g(x-y)f(y)|^p\, dx\,dy\\
\noalign{\vskip6pt}
&\qquad
= \int_{\real^m\times\real^m} |g(y) f(x) - g(-y) f(x+y)|^p\,dx\,dy \\
\noalign{\vskip6pt}
&\qquad
\ge \int_{\real^m} \bigg\{\bigg( \int_{\real^m} |g(y) f(x)-g(-y)f(x+y)|^p\,dx
\bigg)^{1/p}\bigg\}^p\, dy\\
\noalign{\vskip6pt}
&\qquad
\ge \int_{\real^m} \left\{\big|\, |g(y)|\, \|f\|_{L^p(\real^n)}
- |g(-y)|\, \|f\|_{L^p (\real^m)}\big| \right\}^p\,dy \\
\noalign{\vskip6pt}
&\qquad
= \int_{\real^m} \big| \, |g(y)| - |g(-y)|\,\big|^p\,dy
\int_{\real^m} |f(x)|^p\,dx \ .
\end{align*}
By considering $g\ge0$ and the family $\ep^{m/p} f(\ep x)$, one observes
that the inequality is optimal.

On a Lie group $G$ with left-invariant Haar measure $dm$, modular
function $\Delta$ and convolution defined by
\begin{equation*}
(f*g) (w) = \int_G f(z) g(z^{-1} w)\,dm
\end{equation*}
then this lemma becomes:

\begin{NSW-groups}
For $f,g\in L^p (\real^n)$, $1\le p<\infty$ 
\begin{align}
&\int_{G\times G} |g(x^{-1} y) f(x) - g(y^{-1}x) f(y)|^p \,dm\,dm\notag\\
\noalign{\vskip6pt}
&\qquad
\ge \int_G \big|\, |g(y)| - \Delta (y)^{-1/p} |g(y^{-1})|\,\big|^p\,dm
\int_G |f(x)|^p\,dm \ .\label{convolution}
\end{align}
\end{NSW-groups}

\noindent 
One recognizes that the argument is simply the interplay between the 
product structure and the Minkowski inequality for the $L^p$ metric. 
\renewcommand{\qed}{}
\end{proof}

Before beginning the proof of the main theorem, note that
inequality~\eqref{eq1}  reduces to radial functions either by using the
triangle inequality to effect the replacement
\begin{equation*}
f(x) \longrightarrow f_\# (|x|)
= \bigg( \int_{S^{n-1}} |f(|x|\xi)|^p \,d\xi\bigg)^{1/p}
\end{equation*}
or use the symmetrization lemma from \cite{B-PNAS92} to reduce the inequality
%\bibitem{B-PNAS92}		%% 3
to the equimeasurable radial decreasing rearrangement $f^*$, such that
\begin{equation*}
m\{x\in \real^n : f^* (x) >\lambda\}
= m\{x\in \real^n :|f(x)| > \lambda\}\ .
\end{equation*}

\begin{SymLem}
The functional
\begin{equation*}
\Lambda (f,g) = \int_{M\times M} \varphi
\left[ \frac{|f(x) - g(y)|}{\rho (d(x,y))}\right] \kappa [d(x,y)]\,dx\,dy
\end{equation*}
is monotone under equimeasurable radial decreasing rearrangement:
$\Lambda (f,g) \ge \Lambda (f^*,g^*)$.
$M$ is a manifold with distance function $d(x,y)$ and reflection symmetry,
and $\varphi,\rho,\kappa$ are non-negative functions on $[0,\infty)$
with the properties:
{\rm (i)}~$\varphi (0) =0$, $\varphi$ convex and monotone increasing, and
$t\varphi'(t)$ convex;
{\rm (ii)}~$\rho$ monotone  increasing, $\kappa$ monotone decreasing.
\end{SymLem}

\begin{proof}[Proof of Theorem]
Assume that $f$ is now radial; set $t = |x|$, $s = |y|$,
$h(t) = |x|^{\frac{n}p -\beta} f(x)$.
Then inequality~\eqref{eq1} is inequivalent to the following form on
the multiplicative group $\real_+$.
\begin{equation*}
\int_{\real_+\times\real_+} |g(s/t)h(t) - g(t/s)h(s)|^p\ \psi (s/t)\
\frac{ds}s \frac{dt}t
\ge D_{p,\beta} \int_{\real_+} |h(t)|^p\ \frac{dt}t
\end{equation*}
where $g(t) = t^{(n-p\beta)/2p}$,
\begin{equation*}
\psi (t) = \int_{S^{n-1}} \Big[t +\frac1t -25_1\Big]^{-(n+p\beta)/2}\,d\xi
\end{equation*}
and $d\xi$ denotes standard surface measure on $S^{n-1}$.
Note that $\psi$ is symmetric under inversion.
Apply the nonlinear Stein-Weiss lemma using $g\,\psi^{1/p}$ as the second
function in the Lemma, and one finds that
\begin{align*}
D_{p,\beta} & = \int_{\real_+} |t^{\lambda/2} - t^{-\lambda/2}|^p\ \psi (t)
\frac{dt}t\\
\noalign{\vskip6pt}
& = \int_{\real^n} \big| 1- |x|^{-\lambda}\big|^p |x-\eta|^{-n-p\beta}\,dx
\end{align*}
where $\lambda = (n-p\beta)/p$ and $\eta \in S^{n-1}$.
Since the determination of $D_{p,\beta}$ depends only on application of
the Stein-Weiss lemma, the constant must be optimal as observed by a
suitable variation of functions in that inequality.
\renewcommand{\qed}{}
\end{proof}

In contrast to the original Stein-Weiss Lemma (see appendix in 
\cite{B-pitt06}, the key step in the proof of the nonlinear form is 
identified here as the ``lemma''. 
This is largely due to the surprising simplicity and the form of the 
argument plus its natural extension to a non-unimodular Lie groups. 
But for completeness Theorem~\ref{thm8} is reformulated to have the 
``look'' of the original Stein-Weiss lemma with the same proof as above. 

\begin{SW-extended} 
Suppose $K$ is a non-negative symmetric kernel defined on 
$\real^n\times \real^n$, continuous on any domain that excludes the 
diagonal, homogeneous of degree $-n-\gamma$, $K(\delta u,\delta v) = 
\delta^{-n-\gamma} K (u,v)$, $0<\gamma <\min (n,p)$, and 
$K(Ru,Rv) = K(u,v)$ for any $R\in SO(n)$. 
Then for $f\in \S (\real^n)$ and $p\ge 1$, 
\begin{align}%\label{eq:SW-extended}
&\int_{\real^n\times\real^n} |f(x) - f(y)|^p K(x,y)\,dx\,dy 
\ge D_{p,\gamma} \int_{\real^n} |x|^{-\gamma} |f(x)|^p\,dx
\label{eq:SW-extended}\\
\noalign{\vskip6pt}
&\hskip1truein 
D_{p,\gamma} = \int_{\real^n} |1-|x|^{-\lambda} |^p K(x,\eta)\,dx
\notag
\end{align}
for $\lambda = (n-\gamma)/p$ and $\eta\in S^{n-1}$.
This constant is optimal. 
But note that there is no assumption made that it is finite.
\end{SW-extended}

In the context of Maz'ya-Sobolev embedding, Frank and Seiringer 
\cite{FS-arXiv} have given an extension of Theorem~9 to the setting 
of the upper half-space $\real_+^n$ using a similar argument to their 
earlier proof of Theorem~9 \cite{FS}. 
But the nonlinear Stein-Weiss lemma for non-unimodular groups can also be 
used to give a simple proof. 
Note that for both of these inequalities this approach does not give the 
error estimates obtained by Frank and Seiringer. 
To set notation, let $w= (x,y)\in \real_+^n$ with $x\in \real^{n-1}$, $y>0$.

\begin{thm}\label{thm:notation}
For $f\in \S(\real_+^n)$, $0<\beta <1$ and $1\le p < 1/\beta$ 
\begin{align}
&\int_{\real^n_+\times \real_+^n} 
\frac{|f(w) - f(w')|^p}{|w-w'|^{n+p\beta}} \ dw\,dw' 
\ge E_{p,\beta} \int_{\real_+^n} |f(w)|^p \ y^{-p\beta} \,dw
\label{eq:notation}\\
\noalign{\vskip6pt}
&\qquad 
E_{p,\beta} = \frac{\pi^{\frac{n-1}2} \Gamma (\frac{p\beta +1}2)}
{\Gamma (\frac{p\beta +n}2)} 
\int_0^\infty | 1-y^{\beta - \frac1p}|^p \ 
|y-1|^{-1-p\beta}\ dy\ .\notag
\end{align}
This constant is optimal.
\end{thm}

\begin{proof} 
The basic step is to convert this inequality to a problem on 
$n$-dimensional hyperbolic space $\HH^n$ and apply the nonlinear 
Stein-Weiss lemma for groups. 
Here the Liouville-Beltrami upper half-space model is used with: 
Poincar\'e distance
$$d(w,w') = \frac{|w-w'|}{2\sqrt{yy'}\,}$$
left-invariant Haar measure $d\nu = y^{-n}\,dx\,dy$, Riemannian metric 
$ds^2 = y^{-2} (dx^2 + dy^2)$ and invariant gradient $D= y\nabla$.
The group structure of hyperbolic space corresponds to a non-unimodular 
Lie group that is an extension of the affine ``$ax+b$ group.''
Hyperbolic  space $\HH^n$ is identified with the subgroup of 
$SL(n,R)$ given by matrices of the form 
$$\root n\of y 
\begin{pmatrix} \text{\bf 1}&x/y\\ 
\noalign{\vskip6pt}
0&1/y\end{pmatrix}$$
where $x\in \real^{n-1}$ is represented as a column vector and $y>0$. 
Such matrices can act via fractional linear transformation on 
$\real_+^n \simeq \HH^n$
$$w= x+iy\xi \in \real_+^n \longrightarrow \frac{Aw+B}{Cw+D}$$
for a matrix $\left(\begin{smallmatrix} A&B\\ C&D\end{smallmatrix}\right)$
where $A = (n-1) \times (n-1)$ matrix, $B = (n-1) \times 1$ matrix, 
$C = 1 \times (n-1)$ matrix, $D = 1\times 1$ matrix and fixed non-zero 
$\xi\in \real^{n-1}$. 
The group action then corresponds to the multiplication rule 
$$(x,y) (u,v) = (x+yu, yv)$$
for $x,u\in \real^{n-1}$ and $y,v >0$ so this is a non-unimodular group.
Let $\Delta$ denote the modular function defined by 
$\nu (Eg) = \Delta(g) \nu (E)$.
Then $\Delta (x,y) = y^{-(n-1)}$.

Returning to equation~\eqref{eq:notation}, let $F(w) = f(w) y^\lambda$ 
with $\lambda = \frac{n}p -\beta$.
Then this equation becomes 
\begin{gather*}
2^{-n-p\beta} \int_{\HH^n \times\HH^n} 
\frac{|F(w) (y'/y)^{\lambda/2} - F(w') (y/y')^{\lambda/2}|^p}
{d(w,w')^{n+p\beta}}\ d\nu\, d\nu =\\
\noalign{\vskip6pt}
2^{-n-p\beta} \int_{\HH^n\times \HH^n} 
|g(w^{-1} w') F(w) - g(w'{}^{-1}w) F(w')|^p \, d\nu\,d\nu 
\ge E_{p,\beta} \int_{\HH^n} |F|^p\,d\nu
\end{gather*}
where $g(w) = y^{\lambda/2} d(w,\hat 0)^{-(\frac{n}p +\beta)}$ 
with $\hat 0 = (0,1)$. 
Then applying the nonlinear Stein-Weiss lemma for non-unimodular groups 
\begin{align*}
E_{p,\beta} & = \int_{\HH^n} |y^{\lambda/2} - y^{(n-1)/p} y^{-\lambda/2}|^p 
\ y^{(n+p\beta)/2} (x^2 + (y-1)^2)^{-(n+p\beta)/2}\, d\nu\\
\noalign{\vskip6pt}
& = \frac{\pi^{\frac{n-1}2} \Gamma(\frac{p\beta+1}2)}
{\Gamma (\frac{p\beta +n}2)} 
\int_0^\infty |1-y^{\beta -\frac1p}|^p\ 
|y-1|^{-1-p\beta}\ dy\ .
\end{align*}
Note that no information results if the case $p\beta =1$ is treated as 
a limiting  case.
\renewcommand{\qed}{}
\end{proof}

\begin{proof}[Alternate proof]
The interplay between the product structure of the manifold and the 
``triangle inequality'' for $L^p$ norms offers a simpler proof as a direct 
application of Theorem~\ref{thm8}. 
(To better understand this variation in proof, see the discussion on 
page 832 in \cite{B-Fourier97}.) 
\begin{align*}
&\int_{\real_+^n \times\real_+^n} 
\frac{|f(w)- f(w')|^p}{|w-w'|^{n+p\beta}} \,dw\,dw'
\ge \int_{\real_+\times\real_+} |h(y) - h(y')|^p J (y-y')\, dy\,dy'\\
\noalign{\vskip6pt}
&\qquad 
= \int_{\real^{n-1}} \big| 1+ |x|^2\big|^{-n-p\beta}\,dx 
\int_{\real_+\times\real_+} |h(y) - h(y')|^p |y-y'|^{-1-p\beta}\,dy\,dy'\\
\noalign{\vskip6pt}
&\qquad 
\ge \frac{\pi^{\frac{n-1}2} \Gamma (\frac{p\beta +1}2)}
{\Gamma (\frac{p\beta +h}2)} 
\int_0^\infty \left| 1-y^{\beta-\frac1p}\right|^p |y-1|^{-1-p\beta} dy
\int_0^\infty |h(y)|^p y^{-p\beta} dy \\
\noalign{\vskip6pt}
&\qquad 
= \pi^{\frac{n-1}2} \frac{\Gamma(\frac{p\beta+1}2)}
{\Gamma (\frac{p\beta+h}2)} 
\int_0^\infty \left| 1-y^{\beta -\frac1p}\right|^p |y-1|^{-1-p\beta}dy 
\int_{\real_+^n} |f(w)|^p y^{-p\beta} dw
\end{align*}
where 
$$h(y) = \bigg( \int_{\real^{n-1}} |f(x,y)|^p\,dx \bigg)^{1/p}\ ,
\qquad 
J(v) = \int_{\real^{n-1}} |x^2 + v^2|^{-n-p\beta}dx \ .$$

This argument suggests an immediate application to two additional cases: 
the Heisenberg group and manifolds with mixed homogeneity. 
Geometric analysis on the Heisenberg group depends in part on understanding 
the intrinsic $SL(2,R)$ invariance associated with the group manifold. 
But in studying manifolds where homogeneity is ``broken'', it is natural 
to determine if the Besov norm estimates obtained for $\real^n$ and the 
upper half-space $\real_+^n$ can extend simply to the Heisenberg group. 
In the context of Stein-Weiss  integrals, a similar question was 
considered earlier in \cite{B-Fourier97}. 

The Heisenberg group $\H_n$ is realized as the boundary of the Siegel upper
half-space in $\complex^{n+1}$, $D = \{ z\in\complex^{n+1}:\Im z_{n+1} > 
|z_1|^2 +\cdots + |z_n|^2\}$. 
Then $\H_n = \{w = (z,t) : z\in \complex^n$, $t\in\real\}$ with the 
group action 
$$ww' = (z,t) (z',t') = (z+z',\, t+t' +2\Im z\bar z')$$
and Haar measure on the group is given by $dw = dz\,d\bar z\,dt = 4^n\,dx\,dy
\, dt$ where $z= z+iy\in\complex^n$ and $t\in \real$. 
The natural metric here is $d(w,w') = d((z,t),(z',t')) = d(w^{'-1}w,\hat0)$
with 
$$d(w,\hat 0) = d((z,t),0,0)) = \left| |z|^2 +it\right|^{1/2} 
= \left| |z|^4 + t^2\right|^{1/4}\ .$$
\renewcommand{\qed}{}
\end{proof}

\begin{thm}[Besov norms on the Heisenberg group]\label{thm:besov}
Let $f\in \S(\H_n)$, $0<\beta <1$ and $1\le p<2n/\beta$; then 
\begin{gather} 
\int_{\H_n\times\H_n}\frac{|f(w)-f(w')|^p}{d(w,w')^{2n+2+p\beta}} \,dw\,dw'
\ge F_{p,\beta} \int_{\H_n} |z|^{-p\beta} |f|^p\,dw 
\label{eq:besov}\\
\noalign{\vskip6pt}
F_{p,\beta} = \frac{4^n\sqrt{\pi}\, \Gamma [\frac{2n+p\beta}4]}
{\Gamma [\frac{2n+2+p\beta}4]} 
\int_{\real^n} \left| 1-|x|^{-\lambda}\right|^p |x-\eta|^{-2n-p\beta}\,dx 
\notag
\end{gather}
for $\lambda = (2n-p\beta)/p$ and $\eta \in S^{2n-1}$.
\end{thm}

\begin{proof} 
Apply the ``triangle inequality'' to the $t,t'$ integrations
\begin{align*}
&\int_{\H_n\times\H_n} \frac{|f(w)-f(w')|^p}{d(w,w')^{2n+2+p\beta}}\,dw\,dw'\\
\noalign{\vskip6pt}
&\qquad \ge 
\int_{\complex^n\times\complex^n} 
|h(z) - h(z')|^p \int_{\real} \left[ |z-z'|^4 +t^2\right]^{-(2n+2+p\beta)/4}
dt\,dz\,d\bar z\,dz'\, d\bar z'\\
\noalign{\vskip6pt}
&\qquad  = 
\frac{\sqrt{\pi}\, \Gamma [\frac{2n+p\beta}4]}{\Gamma [\frac{2n+2+p\beta}4]}
\int_{\complex^n\times\complex^n} 
\frac{|h(z) - h(z')|^p}{|z-z'|^{2n+p\beta}} 
dz\,d\bar z\,dz'\,d\bar z'\\
\noalign{\vskip6pt}
&\qquad \ge 
\frac{4^n \,\sqrt{\pi}\, \Gamma [\frac{2n+p\beta}4]}
{\Gamma [\frac{2n+2+p\beta}4]} 
\int_{\real^{2n}}\left| 1-|x|^{-\lambda}\right|^p |x-\eta|^{-2n-p\beta} dx
\int_{\complex^n} |z|^{-p\beta} |h|^p \,dx\,d\bar z
\end{align*}
for $\lambda = (2n-p\beta)/p$ and $\eta \in S^{2n-1}$ where 
Theorem~\ref{thm8} is applied to obtain this estimate with 
$$h(z) = \bigg( \int_{\real} |f(z,t)|^p\,dt\bigg)^{1/p}$$
so that inequality~\ref{eq:besov} is obtained. 
\renewcommand{\qed}{}
\end{proof}

A similar result holds for problems with mixed homogeneity (see the 
corresponding discussion for Stein-Weiss potentials on page~1876 in 
\cite{B-pitt06}). 

\begin{thm}[Besov norms with mixed homogeneity]\label{thm:besov-mixed}
Let $f\in \S(\real^{n+m})$, $0<\beta <1$ and $1\le p<m/\beta$, 
$w= (x,v)\in \real^n\times\real^m$; then 
\begin{gather} 
\int_{\real^{n+m}\times\real^{n+m}}\frac{|f(w)-f(w')|}{|w-w'|^{n+m+p\beta}}
dw\,dw' 
\ge G_{p,\beta} \int_{\real^n\times \real^m} |v|^{-p\beta} |f|^p\,dw 
\label{eq:besov-mixed}\\
\noalign{\vskip6pt}
G_{p,\beta} = \frac{\pi^{n/2} \Gamma[\frac{m+p\beta}2]}
{\Gamma [\frac{n+m+p\beta}2]} 
\int_{\real^m} \left| 1-|x|^{-\lambda}\right|^p |x-\eta|^{-m-p\beta}dx
\notag
\end{gather}
for $\lambda = (m-p\beta)/p$ and $\eta \in S^{m-1}$.
\end{thm}

\begin{proof}
Apply the previous argument using the ``triangle inequality'' for the 
$x,x'$ integration and Theorem~\ref{thm8}. 
Observe that $G_{p,\beta}$ is then given by 
$$\int_{\real^n} [1+|x|^2]^{-(n+m+p\beta)/2} dx 
\int_{\real^m} \left| 1-|x|^{-\lambda}\right|^p |x-\eta|^{-m-p\beta} dx$$
with $\lambda,\eta$ as above.
\renewcommand{\qed}{}
\end{proof} 

%%%%%%%%%%%%%%%%%%%%%%%%%%
\mysection{5. Analysis and applications}
\medskip

The Frank-Lieb-Seiringer spectral formula which relates the $L^2$ norm for the 
fractional Laplacian to the weighted Besov norm can be applied in different
ways to give sharp information on embedding questions. 
First, it can be used to obtain global error estimates for Hardy's 
inequality and Pitt's inequality that simultaneously demonstrate sharpness 
and that optimal constants are not attained. 
But more generally, this formula can be combined with the 
Hardy-Littlewood-Sobolev inequality to obtain or recover sharp embedding 
constants. 
Results from \cite{B-ForumMath97} can be viewed as predictive for this 
conceptual framework. 
To illustrate this strategy, the Hardy relation on $\real^3$.
\begin{equation*} 
\int_{\real^3} |\nabla f|^2\,dx 
= \frac14 \int_{\real^3} |x|^{-2} |f|^2\,dx 
+ \int_{\real^3} |\nabla (|x|^{1/2} f)|^2 |x|^{-1} \,dx 
\end{equation*}
can be used to obtain the sharp $H^1(\real)$ Moser inequality calculated   
by Nagy \cite{Nagy} (see also \cite{B-PAnal04})		
\begin{equation}\label{sharpMoser}
\int_{\real} |g|^6 \,dx 
\le \frac4{\pi^2} \int_{\real} |\nabla g|^2\,dx 
\left[\int_{\real} |g|^2\,dx \right]^2\ .
\end{equation}
The Sobolev inequality on $\real^3$ is 
\begin{equation*} 
3 (\pi/2)^{4/3} \bigg(\int_{\real^3}|f|^6\,dx \bigg)^{1/3} 
\le \int_{\real^3} |\nabla f|^2\, dx
\end{equation*}
which results using \eqref{classicalinequality}
\begin{equation*}
3(\pi/2)^{4/3} \bigg( \int_{\real^3} |f|^6 \,dx\bigg)^{1/3}
\le \frac14 \int_{\real^3} |x|^{-2} |f|^2\,dx
+\int_{\real^3} |\nabla (|x|^{1/2} f)|^2 |x|^{-1}\, dx\ .
\end{equation*}
Choose $f$ to be radial, set $t= |x|$ and $h(t) = |x|^{1/2} f(x)$; 
then let $t= e^w$, $g(w) = h(e^w)$ and apply a variational argument to 
obtain \eqref{sharpMoser} 
\begin{equation*}
\int_{\real} |g|^6 \, dw 
\le \frac4{\pi^2} \int_{\real} |\nabla g|^2\, dw 
\bigg[ \int_{\real} |g|^2\,dw \bigg]^2\ .
\end{equation*}
Though the context here is the real line, differentiation is denoted by the 
gradient symbol $\nabla$ to reflect the character of embedding relations.
Using the Hardy-Littlewood-Sobolev inequality, the optimal constant is 
obtained for the extremal function $(\cosh w)^{-1/2}$. 
A different characterization of this approach would be that it provides a 
clever method to discover useful change of variables that develop more 
fully the underlying symmetry of embedding estimates.

This argument can be used with sharp Sobolev embedding on $\real^n$ for 
$n>2$, $p' = 2n/(n-2)$ or $n/p' = n/2\, -\, 1$
\begin{equation*}
\pi n(n-2) \left[ \Gamma (n/2)/\Gamma (n)\right]^{2/n}
\bigg(\int_{\real^n} |f|^{p'} \,dx\bigg)^{2/p'} 
\le \int_{\real^n} |\nabla f|^2\,dx
\end{equation*}
which gives by setting $s=n/p' = n/2\, -\,1$ and using 
\eqref{classicalinequality} 
\begin{gather}
4\pi s(s+1) \left[\Gamma (s+1)/\Gamma (2s+2)\right]^{1/(s+1)} 
\bigg( \int_{\real^n} |f|^{2+\,2/s} \,dx\bigg)^{s/(s+1)} \label{sharpSobolev}\\
\noalign{\vskip6pt}
\le s^2 \int_{\real^n} |x|^{-2} |f|^2\,dx 
+ \int_{\real^n} |\nabla (|x|^2 f)|^2 |x|^{-2s}\,dx\ .\notag
\end{gather}
Choose $f$ to be radial, set $t= |x|$ and $h(t) = |x|^s f(x)$; then 
\begin{gather*}
2^{(2s+1)/(s+1)} s(s+1) \left[\Gamma^2 (s+1)/\Gamma(2s+2)\right]^{1/(s+1)} 
\bigg(\int_0^\infty |h|^{2+\,2/s} \frac{dt}t\bigg)^{s/(s+1)}\\
\noalign{\vskip6pt}
\le s^2 \int_0^\infty |h|^2 \frac{dt}t + \int_0^\infty |th'|^2 \frac{dt}t\ ;
\end{gather*}
now let $t= e^w$ and $g(w) = h(e^w)$ to obtain an $H'$ Sobolev inequality 
on the real line
\begin{gather}
A_s \bigg(\int_\real |g|^{2+\, 2/s} \,dw\bigg)^{s/(s+1)} 
\le s^2 \int_\real |g|^2\, dw + \int_\real |\nabla g|^2\,dw 
\label{HprimeSobolev}\\
\noalign{\vskip6pt}
A_s =2^{(2s+1)/(s+1)} s(s+1)
\left[\Gamma^2 (s+1)/\Gamma (2s+2)\right]^{1/(s+1)}\ .\notag
\end{gather}
Using a standard variational argument $\{g(w) \to g(\delta w)$, $\delta>0\}$, 
one finds the equivalent Gagliardo-Nirenberg inequality which was originally
calculated by Nagy \cite{Nagy}.

\begin{thm}\label{thm5}
For $g\in C^1 (\real)$ and $0<s<\infty$
\begin{gather}
\bigg(\int_\real |g|^{2+\, 2/s}\,dx\bigg)^{2s} 
\le B_s \int_\real |\nabla g|^2\,dx \bigg( \int_\real |g|^2\,dx\bigg)^{2s+1}
\label{eq:thm5}\\
\noalign{\vskip6pt}
B_s = ( s/2)^{2s} (2s+1)^{-(2s+1)} \frac{\Gamma^2 (2s+2)}{\Gamma^4 (s+1)}\ .
\notag
\end{gather}
An extremal function is given by $(\cosh x)^{-s}$.
\end{thm}

\begin{proof}
For the case where $s$ is a half-integer multiple, this inequality  follows 
the previous deduction using sharp Sobolev embedding on $\real^n$. 
The extremal function follows from the analysis on $\real^n$. 
But the general case depends only on observing that an extremal function 
exists by duality for the $H^1$ inequality, and that one can find a unique 
radial-decreasing solution for the Euler-Lagrange variational equation. 
The argument is simple and follows the outline of the Hardy-Littlewood-Sobolev
method used in \cite{B-PAMS01} (see page 1244). 
The two inequalities for $1/p + 1/q =1$, $2<q<q_c = 2n/(n-2)$ for $n>2$, 
$q_c =\infty$ if $n=1,2$ 
\begin{gather}
\bigg( \int_{\real^n} |g|^q\,dx \bigg)^{2/q} 
\le A_q \bigg[\int_{\real^n} |g|^2\,dx + \int_{\real^n} |\nabla g|^2\,dx
\bigg] \label{twoinequalities-1}\\
\noalign{\vskip6pt}
\Big| \int_{\real^n\times\real^n} f(x) G_2 (x-y) f(y)\,dx\,dy\Big| 
\le A_q \left[ \|f\|_{L^p(\real^n)}\right]^2 \label{twoinequalities-2}
\end{gather}
are equivalent by duality and the existence of an extremal function for 
\eqref{twoinequalities-2} will imply existence of an extremal function 
for \eqref{twoinequalities-1}, in particular $g= G_2 * f$ where $G_2$ is 
the Bessel potential defined using the Fourier transform by
\begin{equation*}
\hat G_\alpha (\xi) = (1+4\pi^2 |\xi|^2)^{-\alpha/2}
\end{equation*}
with 
\begin{equation*}
\begin{split}
G_\alpha (x) & = \left[ (4\pi)^{\alpha/2} \Gamma (\alpha/2)\right]^{-1} 
\int_0^\infty e^{-\pi |x|^2/\delta} e^{-\delta /4\pi} 
\delta^{(-n+\alpha)/2} \frac{d\delta}{\delta}\\
\noalign{\vskip6pt}
& = \left[ 2^{(n+\alpha -2)/2} \pi^{n/2} \Gamma (\alpha/2)\right]^{-1} 
|x|^{-(n-\alpha)/2} K_{(n-\alpha)/2} (|x|)\ .
\end{split}
\end{equation*}
Note that $G_\alpha (x)$ is radial decreasing. 
Since $\hat G_\alpha$ is positive, the convolution operator defined by 
$G_\alpha$ is positive-definite and it suffices to study its properties 
on the diagonal. 
By symmetrization, it suffices to treat \eqref{twoinequalities-1} for 
non-negative radial decreasing functions. 
To show the existence of extremal functions, consider a sequence of 
non-negative radial decreasing functions $\{f_m\}$ with $\|f_m\|_p=1$ and
\begin{equation*}
\int_{\real^n\times\real^n} f_m(x) G_2 (x-y) f_m (y)\,dx\,dy 
\xrightarrow[m\to\infty]{} A_q\ .
\end{equation*}
By virtue of the norm condition, $f_m(x) \le c|x|^{-n/p}$. 
Using the Helly selection principle, one can choose a subsequence that 
converges almost everywhere to a function $F$. 
By Fatou's lemma, $\|F\|_p \le 1$. 
But 
\begin{equation*}
f_m (x) G_2 (x-y) f_m (y) \le c^2 |x|^{-n/p} G_2 (x-y) |y|^{-n/p} 
\in L^1 (\real^n \times\real^n)
\end{equation*}
for $2<q<q_c$ where $q$ is the dual exponent to $p$. 
Using the Fourier transform 
\begin{equation*}
\int_{\real^n\times\real^n} |x|^{-n/p}  G_2 (x-y)|y|^{-n/p}\,dx\,dy 
= c\int_{\real^n} |\xi|^{-2n/q} (1+4\pi^2 |\xi|^2)^{-1}\, d\xi <\infty 
\end{equation*}
for $2<q<q_c$. 
This is a nice calculation because it highlights the role of the critical 
index to ensure that the right-hand integral is finite.
Now using the dominated convergence theorem, the existence of an extremal 
function is demonstrated 
\begin{equation*}
\int_{\real^n\times\real^n} f_m (x) G_2 (x-y) f_m(y)\,dx\,dy 
\to \int_{\real^n\times\real^n} F(x) G_2 (x-y) F(y)\,dx\,dy = A_q
\end{equation*}
with $\|F\|_p =1$. 
The existence of an extremal function for \eqref{twoinequalities-1} means 
that an extremal function will exist for the $H^1$ Sobolev embedding estimate 
\eqref{twoinequalities-1}, and in fact one can take $g= G_2 *F$ which is 
bounded and radial decreasing since $G_2 \in L^r (\real^n)$ for $1\le r<q_c$.
Now this extremal function will be a non-negative bounded radial decreasing 
solution of the Euler-Lagrange variational equation
\begin{equation}\label{ELvariation}
-\Delta g +g= c\ g^{q-1}\ .
\end{equation}
For the case $n=1$ where by rescaling the equation becomes 
\begin{equation*}
- g'' + s^2 g = c\ g^{1+\, 2/s}
\end{equation*}
which has a unique bounded symmetric decreasing solution, 
$g(x)=(\cosh x)^{-s}$.
The proof of Theorem~\ref{thm5} is completed by using the extremal function
to compute the constant in \eqref{twoinequalities-2} and then applying 
the variational argument to obtain  \eqref{eq:thm5}. 
This is the only case in which an explicit closed-form solution for the 
extremal function and the constant $A_q$ has been calculated (see Nagy
\cite{Nagy}). 
Weinstein \cite{Weinstein}	%% [28] 
calculated numerically the value of $A_4$ in dimension
two and suggested that his methods could be applied to calculate 
any of the other constants. 
But relatively simple approximations can
be used to obtain good numerical values 
(see \cite{B-PAnal04}, pages 355-357). 
The $n$-dimensional results detailed above are expressed in the 
following theorem.
\renewcommand{\qed}{}
\end{proof}

\begin{thm}\label{ndimensional}
For $f\in L^p(\real^n)$, $g\in C^1 (\real^n)\cap L^2 (\real^n)$ with 
$1/p + 1/q = 1$, $2<q<q_c = 2n/(n-2)$ for $n>2$, $q_c=0\infty$ if $n=1$ or $2$,
the following inequalities are equivalent
\begin{gather*}
\Big| \int_{\real^n\times\real^n} f(x) G_2 (x-y) f(y)\,dx\,dy\Big|
\le A_q  \left[ \|f\|_{L^p (\real^n)} \right]^2\\
\noalign{\vskip6pt} 
\left[\|g\|_{L^q(\real^n)}\right]^2 \le A_q \bigg[ \int_{\real^n} 
|\nabla g|^2 \,dx  + \int_{\real^n} |g|^2\, dx\bigg]\ .
\end{gather*}
Here 
$$\hat G_2 (\xi) = (1+4\pi^2 |\xi|^2)^{-1}\ .$$
Bounded positive radial decreasing extremal functions exist for both 
inequalities. 
For the $H^1$ Sobolev inequality, the extremal will be a radial solution 
of the differential equation 
$$-\Delta g + g =  c\ g^{q-1}\ .$$
\end{thm}

A further interesting remark can be made for the one-dimensional setting 
by using equation~\eqref{HprimeSobolev}. 
Set 
$$d\nu = (\cosh w)^{-2\delta -2}dw\ ,\qquad D= \cosh w\ \frac{d}{dw}$$
and $g= (\cosh w)^{-s} k$. 
Then equation~\eqref{HprimeSobolev} becomes 
\begin{equation}\label{HprimeSobolev-2}
A_s \left[ \|k\|_{L^{2+2/s}(d\nu)}\right]^2 
\le \int |Dk|^2\,d\nu  + s(s+1)\int |k|^2\, d\nu
\end{equation}
which now has some resemblance to Sobolev embedding estimates on a 
curved manifold.

\begin{remarks}
Pitt formulated inequality~\eqref{eq:pitt} for Fourier series.
Zygmund viewed this result as generalizing theorems of Hardy-Littlewood 
and Paley for weighted norms of Fourier coefficients (see notes for 
chapter~12 in {\em Trigonometric Series\/}). 
Extensions of Pitt's inequality for Fourier coefficients of uniformly 
bounded orthonormal systems are given in \cite{S1}, \cite{SW1}. 
Clearly Pitt's inequality encompasses the Hausdorff-Young inequality as well 
as the conformally invariant Hardy-Littlewood-Sobolev inequality. 
A proof for the Euclidean $\real^n$ version of Pitt's inequality 
is given in the Appendix to the author's paper \cite{B-pitt06}. 
Calculation of the best constant for the $L^2$ inequality 
\eqref{eq:pittspectral} was done independently by Herbst \cite{Herbst}, 
Yafaev  \cite{Y} and the author \cite{B-PAMS95}. 
In the recent literature, this inequality has been characterized as a 
Hardy-Rellich inequality or a Hardy-type inequality following the case 
$\alpha=2$. 
The optimal constant for the Stein-Weiss inequality \eqref{eq:steinweiss}
was obtained independently by Samko \cite{Samko} and 
the author \cite{B-pitt06}.
The ``ground state spectral representation'' \eqref{eq:GSR} appears in the 
recent paper by Frank, Lieb and Seiringer \cite{FLS}. 
%But also see the discussion for Lemma~7 by Eilertsen \cite{Eil}. 
A different proof is given in the appendix below. 
The proofs given for the calculation of optimal constants for diagonal 
maps have some overall similarity, but the arguments given by the author 
(\cite{B-PAMS95}, 
\cite{B-Fourier97}, 
\cite{B-weighted06}, 
\cite{B-pitt06}) 
emphasize the geometric symmetry 
corresponding to dilation invariance and characterize the operators in terms
of convolution, as does the proof by Herbst \cite{Herbst}. 
The proofs given in  \cite{KPS}  reflect more the simplicity 
of Schur's  lemma, but an independent argument must be given to show that 
the bound is optimal. 
The method used by Samko \cite{Samko} 
(see Theorem~6.4, page~70 in \cite{KS}) is
to essentially reduce the question to the $n$-dimensional Stein-Weiss 
lemma (see section~2 of Appendix in \cite{B-pitt06}). 
This lemma is a natural extension of Theorem~319 in Hardy, Littlewood 
and P\'olya, {\em Inequalities\/}. 
The advantage of the convolution framework is that sharp constants, 
non-existence of extremals for infinite measures and simplicity of 
iteration are attained directly in one step.
%\begin{schur}
%Let $M$ be a manifold with a $\sigma$-finite measure $d\nu$. 
%Suppose $K$ is a non-negative integral kernel on $M\times M$ with
%$$\int_M K(x,y)\, d\nu (y) 
%= \int_M K(x,y)\, d\nu (x) =A\ .$$
%Then the linear  operator 
%$$(Tf)(x) = \int_M K(x,y) f(y)\,d\nu (y)$$
%is bounded on $L^p(M)$, $1\le p\le\infty$ with 
%$$\|Tf\|_{L^p(M)} \le A\|f\|_{L^p(M)}\ .$$
%For $1<p<\infty$ and a non-finite measure, extremals will not exist.
%\end{schur}
The role of the Besov norm 
\begin{equation}\label{Besovnorm}
C(n,p,\beta) \int_{M\times M} 
\frac{|v(x) -v(y)|^p}{|x-y|^{n+p\beta}}\,dx\,dy\ ,\qquad 
v= (-\Delta)^\alpha f
\end{equation}
with $\dim M=n$, $1\le p<\infty$, $\alpha \ge 0$, $0\le \beta < 1$,  
%\begin{equation*}
%C(n,p,\beta)= \pi^{-n/2\  +  p(1-\beta)} 
%\left[ \frac{\sqrt{\pi} \ \Gamma (\frac{n-p(1-\beta)}2) \Gamma (\frac{n+p}2)}
%{\Gamma (\frac{p+1}2) \Gamma (\frac{n}2) \Gamma(\frac{p(1-\beta)}2)}
%\right]
%\end{equation*}
to characterize Sobolev embedding of fractional order in $L^p(M)$ grew out 
of independent work by Aronszajn, Besov, Calder\'on and Stein (1959--62; 
see especially \cite{Besov61} and \cite{SteinBAMS}), but also 
the significant paper by Gagliardo \cite{Gag}). 
Here the constant $C(n,p,\beta)$ is chosen to facilitate obtaining 
$\int |\nabla v|^p\,dx$ in the limit $\beta\to1$. 
It is interesting that \eqref{Besovnorm} is a characteristic example of a 
larger class of functionals that measure smoothness and are determined by 
the property of monotonicity under equimeasurable radial decreasing 
rearrangement (see Theorem~3 in \cite{B-PNAS92}): 
\begin{equation}\label{radialdecreasing}
\int_{M\times M} \varphi \left[ \frac{|f(x)-g(y)|}{\rho (d(x,y))}\right] 
K [d(x,y)]\,dx\,dy
\end{equation}
where $M$ is a manifold with distance function $d(x,y)$, and $\varphi, \rho,K$
are non-negative functions on $[0,\infty)$ with the properties:
(i)~$\varphi(0)=0$, $\varphi$ convex and monotone increasing; and 
$t\varphi'(t)$ convex;
(ii)~$\rho$ monotone increasing, $K$ monotone decreasing. 
Several sharp examples on $S^n$ comparing the Besov norm with entropy 
are calculated in \cite{B-ForumMath97}. 
An interesting implication from this argument is 
an independent proof of the $L^2$ limit of the Hardy-Littlewood-Sobolev
inequality on the sphere $S^n$.
Asymptotic behavior for Besov norm embedding constants is calculated by both 
Bourgain, Brezis and Mironescu \cite{BBM2000, BBM02} and 
Maz'ya and Shaposhnikova \cite{MS}. 
In extending that work, interesting new $L^p$ Hardy-Rellich inequalities 
with optimal constants and Besov norms controlling fractional 
differentiation have recently been obtained by Frank and Seiringer
\cite{FS, FS-arXiv}.
Stein's ICM lecture at Nice emphasized the importance that analysis on Lie 
groups would play in future development, including the characteristic 
example of $SL(2,R)$ and the fundamental role of dilations, but there was 
not recognition that in going from a manifold to its boundary in the 
noncompact setting, Hadamard manifolds (e.g., spaces with non-positive 
sectional curvature) would  have a central place and explicit calculations 
would need estimates for non-unimodular groups. 
\end{remarks}

%%%%%%%%%%%%%%%%%%%%%%%%%%%%%
\section*{Appendix}

To make the present discussion more complete, quick calculations are 
provided to obtain the representation formulas for the fractional Laplacian.

\begin{classical}[Aronszajn-Smith]
For $f\in \S(\real^n)$, $0<\alpha <2$ 
\begin{equation}\label{eq:classical}
\int_{\real^n\times\real^n} \frac{|f(x)-f(y)|^2}{|x-y|^{n+\alpha}}
\, dx\,dy 
= D_\alpha \int_{\real^n} |\xi|^\alpha |\hat f(\xi)|^2\,d\xi
\end{equation}
\begin{equation*}
D_\alpha = \frac4{\alpha} \pi^{\frac{n}2 + \alpha}
\frac{\Gamma (1-\frac{\alpha}2)}{\Gamma (\frac{n+\alpha}2)}
\end{equation*}
\end{classical}

\begin{proof} 
This is a simple application of the Plancherel theorem (see Stein, page 140).
\begin{equation*}
\begin{split}
\int_{\real^n\times\real^n} 
\frac{|f(x) - f(y)|^2}{|x-y|^{n+\alpha}} \,dx\,dy 
& = \int_{\real^n} \frac1{|w|^{n+\alpha}} 
\bigg[ \int_{\real^n} |f(x+w) - f(x)|^2\,dx\bigg]\,dw\\
\noalign{\vskip6pt}
& = \int_{\real^n} \frac1{|w|^{n+\alpha}} \int_{\real^n} 
|e^{2\pi iw\cdot\xi} - 1|^2\ 
|\hat f(\xi)|^2\,dx\,dw\\
\noalign{\vskip6pt}
& = \int_{\real^n} \frac1{|w|^{n+\alpha}} |e^{2\pi iw\cdot\eta} -1|^2\,dw 
\int_{\real^n} |\xi|^\alpha |\hat f(\xi)|^2\,d\xi
\end{split}
\end{equation*}
with $\eta \in S^{n-1}$. 
Then 
\begin{equation*}
\begin{split}
\int_{\real^n} \frac1{|w|^{n+\alpha}} |e^{2\pi iw\cdot\eta} -1|^2\,dw 
& = 2\int_{\real^n} \frac1{|w|^{n+\alpha}} (1-\cos 2\pi w\cdot\eta)\,dw\\
\noalign{\vskip6pt}
& = \frac{2\pi^{\frac{n+\alpha}2}}{\Gamma(\frac{n+\alpha}2)} 
\int_{\real^n} (1-\cos 2\pi w\cdot\eta) 
\int_0^\infty t^{\frac{n+\alpha}2 -1} e^{-\pi +w^2}\,dt \\
\noalign{\vskip6pt}
& = \frac{2\pi^{\frac{n+\alpha}2}}{\Gamma (\frac{n+\alpha}2)}
\int_0^\infty t^{\frac{n+\alpha}2 -1} \int_{\real^n} 
(1-\cos 2\pi w\cdot\eta) e^{-\pi tw^2}\,dw\\
\noalign{\vskip6pt}
& = \frac{2\pi^{\frac{n+\alpha}2}}{\Gamma(\frac{n+\alpha}2)} 
\int_0^\infty \mkern-12mu 
t^{\frac{\alpha}2 -1} (1-e^{-\pi/t})\,dt 
= \frac{2\pi^{\frac{n+\alpha}2}}{\Gamma(\frac{n+\alpha}2)} 
\int_0^\infty \mkern-12mu t^{-\frac{\alpha}2 -1} (1-e^{-t})\,dt\\
\noalign{\vskip6pt}
& = \frac4{\alpha} 
\frac{\pi^{\frac{n}2 +\alpha}}{\Gamma(\frac{n+\alpha}2)} 
\int_0^\infty t^{-\alpha/2} e^{-t} \,dt 
= \frac4{\alpha} \pi^{\frac{n}2 +\alpha} 
\frac{\Gamma (1-\frac{\alpha}2)}{\Gamma(\frac{n+\alpha}2)}\ .
\end{split}
\end{equation*}
The positivity of the integrands justify the exchange of orders of 
integration using Fubini's theorem. 
An alternative argument can be given using distribution theory and 
Green's theorem.
\begin{equation*}
\begin{split}
&2\int_{\real^n} \frac1{|w|^{n+\alpha}} (1-\cos 2\pi w\cdot\eta)\,dw\\
\noalign{\vskip6pt}
&\qquad = 
\left[ \alpha \Big(\frac{n+\alpha}2 -1\Big)\right]^{-1} 
\int_{\real^n} \Delta \left(\frac1{|w|^{n+\alpha -2}}\right) 
(1-\cos 2\pi w\cdot \eta)\,dw\\
\noalign{\vskip6pt}
&\qquad = 
\left[ \alpha\Big( \frac{n+\alpha}2 -1\Big)\right]^{-1} 
\int_{\real^n} \frac1{|w|^{n+\alpha -2}} \Delta (1-\cos 2\pi w\cdot\eta)\,dw\\
\noalign{\vskip6pt}
&\qquad = 
4\pi^2 \left[\alpha \Big(\frac{n+\alpha}2-1\Big)\right]^{-1} 
\int_{\real^n} \frac1{|w|^{n+\alpha -2}} \cos 2\pi w\cdot \eta\, dw\\
\noalign{\vskip6pt}
&\qquad = 
4\pi^2 \left[\alpha \Big(\frac{n+\alpha}2 -1\Big)\right]^{-1} 
{\mathcal F} \Big[ \frac1{|w|^{n+\alpha -2}}\Big] (\eta)\\
\noalign{\vskip6pt}
&\qquad = 
\frac{4\pi^{\frac{n}2 +\alpha}}{\alpha}\ \  
\frac{\Gamma (1-\frac{\alpha}2)}{\Gamma (\frac{n+\alpha}2)}
\end{split}
\end{equation*}
\end{proof}

\begin{weighted}[Frank-Lieb-Seiringer]
For $f\in \S(\real^n)$,  $0<\alpha <\min (2,n)$ and 
$g(x) = |x|^\lambda f(x)$, $0<\lambda < n-\alpha$ 
\begin{equation}\label{eq:weighted}
\begin{split}
D_\alpha \int_{\real^n} |\xi|^\alpha |\hat f(\xi)|^2\, d\xi 
&= \int_{\real^n\times\real^n} 
\frac{|g(x)-g(y)|^2}{|x-y|^{n+\alpha}} \ |x|^{-\lambda} |y|^{-\lambda}
\, dx\,dy\\
\noalign{\vskip6pt}
&\qquad + \Lambda (\alpha,\lambda,n) \int_{\real^n} |x|^{-\alpha} 
|f(x)|^2\,dx\\
\noalign{\vskip6pt}
\Lambda (\alpha,\lambda,n) &= \pi^{-\alpha}\ D_\alpha 
\left[ \frac{\Gamma (\frac{n-\lambda}2)\Gamma(\frac{\lambda+\alpha}2)}
{\Gamma(\frac{\lambda}2)\Gamma(\frac{n-\alpha-\lambda}2)}
\right]
\end{split}
\end{equation}
\begin{equation*}
D_\alpha = \frac4{\alpha} \ \pi^{\frac{n}2 +\alpha} \ 
\frac{\Gamma (1-\frac{\alpha}2)}{\Gamma (\frac{n+\alpha}2)}
\end{equation*}
For $\lambda = \frac{n-\alpha}2$, then $\Lambda (\alpha,\lambda,n) = 
D_\alpha/C_\alpha$
$$C_\alpha = \pi^\alpha \left[ \frac{\Gamma (\frac{n-\alpha}4)}
{\Gamma(\frac{n+\alpha}4)}\right]^2$$
\end{weighted}

\begin{proof}
Using the classical formula
\begin{equation*}
\begin{split}
&D_\alpha \int_{\real^n} |\xi|^\alpha |\hat f(\xi)|^2\,d\xi 
= \int_{\real^n\times\real^n}
\frac{|f(x) - f(y)|^2}{|x-y|^{n+\alpha}}\, dx\,dy\\
\noalign{\vskip6pt}
&\qquad
= \int_{\real^n\times\real^n} 
\frac{|g(x) -g (y)|^2}{|x-y|^{n+\alpha}} \ |x|^{-\lambda} |y|^{-\lambda}
\,dx\,dy 
+ 2 \int_{\real^n\times\real^n} \left[ 1-\frac{|x|^\lambda}{|y|^\lambda}
\right] 
\frac{|f(x)|^2}{|x-y|^{n+\alpha}} \,dx\,dy\\
\noalign{\vskip6pt}
&\qquad 
= \int_{\real^n\times\real^n} 
\frac{|g(x) - g(y)|^2}{|x-y|^{n+\alpha}}\  
|x|^{-\lambda} |y|^{-\lambda}\, dx\, dy 
+ \Lambda (\alpha,\lambda,n) \int_{\real^n} |x|^{-\alpha} |f(x)|^2\,dx
\end{split}
\end{equation*}
$$\Lambda (\alpha,\lambda,n) = 2 \int_{\real^n} 
\left( 1-\frac1{|y|^\lambda}\right) \frac1{|y-\eta|^{n+\alpha}} \, dy\ ,
\qquad \eta \in S^{n-1}$$
Note that $\int_{S^{n-1}} (1-|y+n|^{-\lambda})\,d\eta \simeq O(|y|^2)$ 
as $|y|\to0$ so $\Lambda (\alpha,\lambda,n)$ is well-defined for 
${0<\alpha <2}$.
\end{proof}

\begin{lem} 
For $0< \alpha < \min (2,n)$, $0<\lambda <n-\alpha$ and $\eta\in S^{n-1}$
\begin{equation}\label{eq34}
2\int_{\real^n} \left( 1-\frac1{|y|^\lambda}\right) 
\frac1{|y-\eta|^{n+\alpha}}\, dy 
= \pi^{-\alpha} D_\alpha
\left[ \frac{\Gamma (\frac{n-\lambda}2) \Gamma (\frac{\lambda+\alpha}2)}
{\Gamma (\frac{\lambda}2) \Gamma (\frac{n-\alpha-\lambda}2)}\right]\ .
\end{equation}
\end{lem}

\begin{proof}[First Proof]
For $n\ge3$, this constant can be calculated using the method of 
distribution theory, Green's theorem and analytic continuation. 
Initially let $0<\lambda < n-2$

\begin{equation*}
\begin{split}
&2\int_{\real^n} \left( 1-\frac1{|y|^\lambda}\right) 
\frac1{|y-\eta|^{n+\alpha}}\, dy 
= \frac2{\alpha (n+\alpha-2)} \int_{\real^n}
\left( 1-\frac1{|y|^\lambda}\right) \Delta 
(|y-\eta|^{-n-\alpha +2})\,dy\\
\noalign{\vskip6pt}
&\qquad 
= \left[\alpha \Big(\frac{n+\alpha}2 -1\Big)\right]^{-1}
\int \Delta \left( 1-\frac1{|y|^\lambda}\right) 
|y-\eta|^{-n-\alpha +2}\, dy\\
\noalign{\vskip6pt}
&\qquad 
= \lambda (n-2-\lambda) \left[\alpha\Big( \frac{n+\alpha}2 -1\Big)\right]^{-1}
\int_{\real^n} |y|^{-\lambda-2} |y-\eta|^{-n-\alpha +2}\,dy\\
\noalign{\vskip6pt}
&\qquad 
=\lambda (n-2-\lambda) \left[\alpha \Big(\frac{n+\alpha}2 -1\Big)\right]^{-1}
(|y|^{-\lambda-2} * |y|^{-n-\alpha +2}) (\eta)\\
\noalign{\vskip6pt}
&\qquad 
= \pi^{n/2} \lambda (n-2-\lambda) 
\left[\alpha \Big(\frac{n+\alpha}2 -1\Big)\right]^{-1}
\left[ \frac{\Gamma (\frac{n-2-\lambda}2) \Gamma(\frac{2-\alpha}2)
\Gamma (\frac{\lambda +\alpha}2)}
{\Gamma(\frac{\lambda+2}2) \Gamma(\frac{n+\alpha-2}2)
\Gamma (\frac{n-\alpha-\lambda}2)} \right]\\
\noalign{\vskip6pt}
&\qquad 
= \frac{4\pi^{n/2}}{\alpha} 
\left[\frac{\Gamma(\frac{n-\lambda}2) \Gamma (1-\frac{\alpha}2) 
\Gamma (\frac{\lambda+\alpha}2)} 
{\Gamma(\frac{\lambda}2) \Gamma(\frac{n+\alpha}2) 
\Gamma (\frac{n-\alpha-\lambda}2)} \right] \\
\noalign{\vskip6pt}
&\qquad 
= \pi^{-\alpha} D_\alpha 
\left[\frac{\Gamma (\frac{n-\lambda}2) \Gamma (\frac{\lambda+\alpha}2)}
{\Gamma (\frac{\lambda}2) \Gamma (\frac{n-\alpha -\lambda}2)}\right]\ .
\end{split}
\end{equation*}
Since both the beginning and final terms are analytic in $\lambda$ on 
the strip $0< Re\, \lambda <n-\alpha$, this formula extends by analytic 
continuation to hold on that region for $n\ge 3$:
$$2\int_{\real^n} \left( 1-\frac1{|y|^\lambda}\right) 
\frac1{|y-\eta|^{n+\alpha}} \,dy 
= \pi^{-\alpha} D_\alpha 
\left[ \frac{\Gamma(\frac{n-\lambda}2) \Gamma(\frac{\lambda+\alpha}2)}
{\Gamma(\frac{\lambda}2) \Gamma(\frac{n-\alpha-\lambda}2)}\right]\ .$$
Set $\lambda = \sigma +\frac{n-\alpha}2$; then for $-(\frac{n-\alpha}2) 
<\sigma < \frac{n-\alpha}2$
$$\Lambda \left(\alpha,\sigma +\frac{n-\alpha}2,\eta\right) 
= \pi^{-\alpha} D_\alpha 
\left[\frac{\Gamma (\frac{n+\alpha}4 - \frac{\sigma}2) 
\Gamma (\frac{n+\alpha}4 +\frac{\sigma}2)} 
{\Gamma (\frac{n-\alpha}4 +\frac{\sigma}2) 
\Gamma (\frac{n-\alpha}4 - \frac{\sigma}2)} \right]\ .$$
Observe that $\Lambda (\alpha,\sigma +\frac{n-\alpha}2,n)$ as a function 
of $\sigma$ is symmetric and has a negative second derivative so $\Lambda$ 
is concave in $\sigma$ for the allowed range and has a maximum at 
$\sigma=0$ or $\lambda = (n-\alpha)/2$.
Hence 
$$\Lambda (\alpha ,\lambda,n) \le \pi^{-\alpha} D_\alpha 
\left[ \frac{\Gamma (\frac{n+\alpha}4)}{\Gamma(\frac{n-\alpha}4)}\right]^2 
= D_\alpha /C_\alpha\ .$$
\end{proof}

\begin{proof}[Second Proof] 
To give a full proof of the lemma, an integral representation is used 
for the factor $1-|y|^{-\lambda}$:
$$1-|y|^{-\lambda} = \frac{\pi^{\lambda/2}}{\Gamma (\lambda/2)} 
\int_0^\infty t^{\frac{\lambda}2 -1} 
(e^{-\pi t} - e^{-\pi t|y|^2})\,dt\ .$$
Then
\begin{equation*}
\begin{split}
&2\int_{\real^n} (1-|y|^{-\lambda}) |y-\eta|^{-n-\alpha} \,dy\\
\noalign{\vskip6pt}
&\qquad 
= \frac{2\pi^{\lambda/2}}{\Gamma (\lambda/2)} \int_{\real^n}
\int_0^\infty t^{\frac{\lambda}2 -1} 
(e^{-\pi t} -e^{-\pi t|y|^2}) \,dt\ |y-\eta|^{-n-\alpha}\,dy\\
\noalign{\vskip6pt}
&\qquad 
= \frac{2\pi^{\lambda/2}}{\Gamma (\lambda/2)} \int_0^\infty 
t^{\frac{\lambda}2 -1} \int_{\real^n} 
(e^{-\pi t} - e^{-\pi t|y|^2}) |y-\eta|^{-n-\alpha}\,dy\,dt\\
\noalign{\vskip6pt}
&\qquad 
= \frac{2\pi^{\lambda/2}}{\alpha (n+\alpha-2)\Gamma(\frac{\lambda}2)}
\int_0^\infty t^{\frac{\lambda}2-1} \int_{\real^n} 
(e^{-\pi t} - e^{-\pi t|y|^2}) \Delta (|y-\eta|^{-n-\alpha+2})\,dy\,dt \\
\noalign{\vskip6pt}
&\qquad 
= \frac{2^{\pi^{\lambda/2}}}{\alpha(n+\alpha-2)\Gamma (\frac{\lambda}2)}
\int_0^\infty t^{\frac{\lambda}2-1} \int_{\real^n} \Delta 
(-e^{-\pi t|y|^2}) |y-\eta|^{-n-\alpha+2}\,dy\,dt \\
\noalign{\vskip6pt}
&\qquad 
= \frac{4\pi^{\frac{n}2 +\frac{\lambda}2 +\alpha}\ \Gamma (\frac{1-\alpha}2)}
{\alpha \Gamma (\frac{\lambda}2) \Gamma (\frac{n+\alpha}2)}
\int_0^\infty t^{\frac{\lambda}2 - \frac{n}2 -1} \int_{\real^n}
e^{2\pi i\xi\cdot\eta} e^{-\pi |\xi|^2/t} |\xi|^\alpha\,d\xi\, dt\\
\noalign{\vskip6pt}
&\qquad 
= D_\alpha \frac{\pi^{\lambda/2}}{\Gamma(\frac{\lambda}2)} \int_{\real^n}
e^{2\pi i\xi\cdot\eta} \int_0^\infty t^{\frac{\lambda}2 -\frac{n}2-1}  
e^{-\pi |\xi|^2/t} \,dt\, |\xi|^\alpha\,d\xi\\
%\noalign{\vskip6pt}
\end{split}
\end{equation*}
\begin{equation*}
\begin{split}
&\qquad 
= D_\alpha \frac{\pi^{\lambda/2}}{\Gamma(\frac{\lambda}2)} \int_{\real^n}
e^{2\pi i\xi\cdot\eta} \int_0^\infty t^{\frac{n-\lambda}2 -1} 
e^{-\pi t|\xi|^2} \,dt \, |\xi|^\alpha\,d\xi
\hskip1truein\\
\noalign{\vskip6pt}
&\qquad 
= D_\alpha \pi^{-\frac{n}2 +\lambda}\  
\frac{\Gamma (\frac{n-\lambda}2)}{\Gamma(\frac{\lambda}2)} \int_{\real^n}
e^{2\pi i\xi\cdot\eta} \frac1{|\xi|^{n-\alpha-\lambda}}\,d\xi\\
\noalign{\vskip6pt}
&\qquad 
= D_\alpha \pi^{-\alpha} \ \frac{\Gamma(\frac{\alpha+\lambda}2)}
{\Gamma (\frac{n-\alpha-\lambda}2)}\ 
\frac{\Gamma(\frac{n-\lambda}2)}{\Gamma(\frac{\lambda}2)} 
= \Lambda (\alpha,\lambda,n)\ .
\end{split}
\end{equation*}
The last exchange of orders of integration requires that the calculation 
be done in the context of distributions.
Different proofs are given in \cite{FLS}.
\end{proof}

%%%%%%%%%%%%%%%%%%%%%%%%%%%%%
In terms of the modulus of continuity, the Aronszajn-Smith formula can be 
extended to larger values of $\alpha$ (see Stein, pages 140, 162--163).

\begin{classical}[Stein] 
For $f\in \S(\real^n)$, $0<\alpha <4$
\begin{gather}
\int_{\real^n\times\real^n} 
\frac{|f(x+y) + f(x-y) -2f(x)|^2}{|y|^{n+\alpha}} dx\,dy 
= E_\alpha \int_{\real^n} |\xi|^\alpha |\hat f(\xi)|^2\,d\xi 
\label{eq:stein}\\
\noalign{\vskip6pt}
E_\alpha  = \frac{4-2^\alpha}{2-\alpha}\ \frac8{\alpha}\ 
\pi^{\frac{n}2 +\alpha} \ 
\frac{\Gamma (2-\frac{\alpha}2)}{\Gamma (\frac{n+\alpha}2)}\ .
\notag
\end{gather}
\end{classical}

\begin{proof} 
Apply the Plancheral formula
\begin{align*}
& \int_{\real^n\times\real^n} 
\frac{|f(x+y) + f(x-y)-2f(x)|^2}{|y|^{n+\alpha}}\,dx\,dy\\
\noalign{\vskip6pt}
&\qquad =  
\int_{\real^n} \frac1{|y|^{n+\alpha}} \int_{\real^n} 
\left| e^{2\pi iy\cdot\xi} + e^{-2\pi iy\cdot\xi} -2\right|^2 
|\hat f(\xi)|^2\,d\xi\,dy\\
\noalign{\vskip6pt} 
&\qquad = 
\int_{\real^n} |\hat f(\xi)|^2 \int_{\real^n} \left| e^{2\pi iy\cdot\xi} 
+ e^{-2\pi i y\cdot \xi} -2\right|^2 \ \frac1{|y|^{n+\alpha}}\,dy\,d\xi\\
\noalign{\vskip6pt}
&\qquad = 
\int_{\real^n} \frac1{|y|^{n+\alpha}} 
\left| e^{2\pi iy\cdot\eta} + e^{-2\pi iy\cdot\eta} -2\right|^2\,dy 
\int_{\real^n} |\xi|^\alpha |\hat f(\xi)|^2\,d\xi
\end{align*}
for $\eta \in S^{n-1}$. 
Observe that 
\begin{gather*}
\left| e^{2\pi iy\cdot \eta} + e^{-2\pi iy\cdot\eta} -2\right|^2 
= 4 | 1-\cos 2\pi y\cdot \eta|^2\\
\noalign{\vskip6pt}
= 8 (1-\cos 2\pi y\cdot\eta) - 2(1-\cos 4\pi y\cdot\eta)\ .
\end{gather*}
Using the previous calculation of $D_\alpha$, $E_\alpha = (4-2^\alpha)D_\alpha$
and 
$$E_\alpha = \frac{4-2^\alpha}{2-\alpha}\ \frac8{\alpha}\ 
\pi^{\frac{n}2 +\alpha}\ 
\frac{\Gamma (2-\frac{\alpha}2)}{\Gamma (\frac{n+\alpha}2)}\ .$$

This result leads to three interesting formulas:
\begin{equation}
\text{(i)}\hskip.20truein
\int_{\real^n\times\real^n} 
\frac{|f(x+y) + f(x-y)-2f(x)|^2}{|y|^{n+2}}\,dx\,dy 
= \frac{\ln 2}{n}\ \frac{2\pi^{n/2}}{\Gamma (n/2)} \ \int_{\real^n} 
|\nabla f|^2\,dx 
\hskip.60truein
\label{int1}
\end{equation}
\begin{equation}
\text{(ii)}\hskip.25truein
\int_{\real^n\times\real^n}
\frac{|f(x+y) + f(x-y) - 2f(x)|^2}{|y|^{n+\alpha}}\,dx\,dy 
\ge \frac{E_\alpha}{C_\alpha} \ \int_{\real^n} |x|^{-\alpha} 
|f(x)|^2\,dx 
\hskip.65truein
\label{int2} 
\end{equation} 
for $0<\alpha < \min (4,n)$ and 
$$C_\alpha = \pi^\alpha \left[ \Gamma \Big(\frac{n-\alpha}4\Big) \Big\slash 
\Gamma \Big(\frac{n+\alpha}4\Big)\right]^2\ ,$$
and 
\begin{equation*}%\label{int3}
\text{(iii)}\hskip.15truein
\int_{\real^n\times\real^n} \mkern-36mu 
\frac{|f(x\!+\!y)\! +\! f(x\!-\!y)\! -\! 2f(x)|^2}{|y|^{n+\alpha}}\,dx\,dy 
\ge \int_{\real^n\times\real^n}\mkern-36mu 
\frac{|f^* (x\!+\!y)\! +\! f^* (x\!-\!y)\! -\! 2f^* 
(x)|^2}{|y|^{n+\alpha}}\,dx\,dy 
\hskip.20truein (48)
\end{equation*}
for $0<\alpha \le 2$ and $f^*$ is the equimeasurable radial decreasing 
rearrangement of $|f|$ on $\real^n$.
\renewcommand{\qed}{}
\end{proof}

%%%%%%%%%%%%%%%%%%%%%%%%%%%%%
\section*{Acknowledgement}

I would like to thank Michael Perelmuter for drawing my attention to the 
arguments of Herbst \cite{Herbst} and Kovalenko-Perelmuter-Sememov 
\cite{KPS}, Sharif Nasibov for reference~\cite{Nas},
 and Emanuel Carneiro for helpful comments.
I am very indebted to Eli Stein for the many ideas he shared in his 
marvelous lectures on Fourier Analysis in the fall of 1969 in new Fine Hall 
and in his mimeographed notes which were to become his classic text on 
{\em Singular Integrals\/}. 
Some of these ideas return in this paper. 

%%%%%%%%%%%%%%%%%%%%%%%%%%%%%

\end{document}